\documentclass{article}
\usepackage[utf8]{inputenc}
\usepackage{amsmath,amssymb,amsthm,todonotes}
\usepackage{amsfonts}
\usepackage{booktabs}
\usepackage{graphics} 
\usepackage{siunitx}
\usepackage{changes}

\definecolor{ForestGreen} {cmyk}{0.91,0,0.88,0.12}
\definechangesauthor[color=magenta]{SB}
\definechangesauthor[color=red]{PA}
\definechangesauthor[color=ForestGreen]{DP}
\usepackage{multirow}
\usepackage{subfig}
\usepackage{tikz}
\usepackage{pgfplots}
\usepackage{braket}
\usepackage{authblk}

\pgfplotsset{compat=newest}
\usepackage{hyperref}

\hypersetup{
    colorlinks=true,
    linkcolor=blue,
    filecolor=magenta,      
    urlcolor=cyan,
}
\usepackage{comment}
\tikzset{immagine/.style={%
  above right, inner sep=0pt, outer sep=0pt}}

\newcommand{\roundPrecision}{2}

\theoremstyle{plain}
\newtheorem{theorem}{Theorem}[section]

\newtheorem{corollary}[theorem]{Corollary}

\theoremstyle{definition}

\newtheorem{remark}[theorem]{Remark}

\theoremstyle{plain}

\newcommand{\argmin}{\mathop{\mathrm{arg\,min}}}

\newcommand{\Tau}{\mathcal{T}}

\newcommand{\E}  {\mathcal{E}}

\newcommand{\V}{V}

\newcommand{\EOD}{\end{document}}

\begin{document}

\title{The Virtual Element Method \\ for a Minimal Surface Problem}
 
  \author[1]{P. F. Antonietti}
  \author[2]  {S. Bertoluzza}
    \author[3]{D. Prada}
  \author[4]{M. Verani}
  
\affil[1]{
    MOX, Dipartimento di Matematica,
    Politecnico di Milano, Italy;
    \emph{e-mail: paola.antonietti@polimi.it}
  }
  \affil[2]{Istituto di Matematica Applicata e Tecnologie Informatiche - CNR, Pavia, Italy
  \emph{e-mail: silvia.bertoluzza@imati.cnr.it}
  }
    \affil[3]{Istituto di Matematica Applicata e Tecnologie Informatiche - CNR, Pavia, Italy
    \emph{e-mail: daniele.prada@imati.cnr.it}
  }
  \affil[4]{
    MOX, Dipartimento di Matematica,
    Politecnico di Milano, Italy and 
    Istituto di Matematica Applicata e Tecnologie Informatiche - CNR, Pavia, Italy
    
    \emph{e-mail: marco.verani@polimi.it}
  }
\maketitle
\begin{abstract}
In this paper we consider the {V}irtual {E}lement discretization of a minimal surface problem, a quasi-linear elliptic partial differential equation modeling the problem of minimizing the area of a surface subject to a prescribed boundary condition. We derive optimal error estimate and present several numerical tests assessing the validity of the theoretical results. 
\end{abstract}

\section{Introduction}

In recent years, the numerical approximation of partial differential equations on computational meshes composed by arbitrarily-shaped
polygonal/polyhedral (polytopal, for short) elements has been 
the subject of an intense research activity.
Examples of such methods include
the Mimetic Finite Difference method,
the Polygonal Finite Element Method,
the polygonal Discontinuous Galerkin Finite Element Methods,
the Hybridizable Discontinuous Galerkin and Hybrid High-Order Methods,
the Gradient Discretization method,
the Finite Volume Method,
the BEM-based FEM,
 the Weak Galerkin method 
 and the {V}irtual {E}lement method (VEM). For more details see the special issue \cite{special-issue} and the references therein. VEM has been  introduced in \cite{VEMvolley} for elliptic problems and later extended to several different linear and non-linear differential problems. While the analysis of linear problems is much more flourished, the study of {V}irtual {E}lement discretization for non-linear problems is much less developed (see, e.g., \cite{absv_VEM_cahnhilliard,Gatica-1,Gatica-2,Vacca:2018b,Lovadina-1,Lovadina-2,Cangiani-quasilinear,obstacle,Adak:1,Adak:2,Adak:3,Liu:2019}).
 In this paper we contribute to fill this gap by addressing the (lowest order) {V}irtual {E}lement discretization of a minimal surface problem (see, e.g., \cite{Ciarlet:book} for its finite element discretization).
 More precisely, in Section~\ref{Sec:1} we introduce the continuous problem together with its {V}irtual {E}lement discretization, while in Section~\ref{Sec:2} we derive optimal error estimate in the  $H^1$-norm, {under a condition on the discrete solution, the validity of which can be checked ``a posteriori''}. Finally, in Section~\ref{Sec:3} we present several numerical results assessing the validity of the theoretical estimate {and confirming that optimal convergence is indeed achieved}. Moreover, the convergence properties in the $L^2$-norm is numerically investigated.

\subsection{Notation}
Throughout the paper we shall use the standard notation of the Sobolev spaces $H^m(\mathcal{D})$ for 
a nonnegative integer $m$ and an open bounded domain $\mathcal{D}$. The $m$-th seminorm 
of the function $v$ will be denoted by
$$ \vert v \vert{^2}_{m,\mathcal{D}} = \sum_{\vert \alpha \vert = m} 
\bigg\| \frac{\partial^{\vert \alpha\vert} v}{\partial_{x_1}^{\alpha_1} \partial_{x_2}^{\alpha_2}}\bigg\|^2_{0,\mathcal{D}},$$
where 
$\| \cdot \|_{0,\mathcal{D}}$ stands for the $L^2(\mathcal{D})$ norm and we set
$\vert {\alpha}\vert = \alpha_1 +\alpha_2$ 
for the nonnegative multi-index ${\alpha} =(\alpha_1,\alpha_2)$.
For any integer $m\geq 0$, $\mathbb{P}^m ({\mathcal D})$ is the space of polynomials of total degree up to $m$ defined on ${\mathcal{D}}$.
Moreover, $n=(n_1,n_2)$ is the outward unit normal vector to $\partial \mathcal{D}$, the boundary of $\mathcal{D}$. Finally,
we will employ the symbol $\lesssim$ for an inequality 
holding up to a constant independent of the mesh size.

\section{Continuous problem and {its} VEM discretization}\label{Sec:1}
Let $\Omega\subseteq \mathbb{R}^2$ be a bounded open set. 
In the following,  we will employ the following notation 
$$f(\cdot)=\sqrt{1+\vert \nabla(\cdot)\vert^2}.$$ 
Let $\varphi$ be a function given on the boundary $\Gamma=\partial\Omega$. The minimal surface problem amounts to finding a function $u$ which minimizes the functional  $$J(v)=\int_{\Omega} f(v)dx$$ over a suitable space of functions which are equal to $\varphi$ on $\Gamma$. The existence and uniqueness of a solution is a delicate mathematical issue (see, e.g., \cite{Ciarlet:book} and the references therein). Here, with the aim of simplifying the analysis, we follow the framework considered, e.g., in \cite{Ciarlet:book} and make the following hypotheses: the domain $\Omega$ is a convex polygonal set and the function $\varphi$ is the trace over $\Gamma$ of a function (by abuse of notation still denoted by $\varphi$) of $H^2(\Omega)$. Moreover, for the subsequent discussion, as in \cite{Ciarlet:book}, we consider that the minimal surface problem consists in solving the following:

\begin{equation}\label{min}
u=\argmin_{v\in V^\varphi}J(v),
\end{equation}
where  $V^\varphi=\{v\in H^1(\Omega): ~ u=\varphi \text{~on~} \partial\Omega\}$.
Note that $u$ is the solution to \eqref{min} if and only if 
$u\in V^\varphi$ solves
\begin{equation}\label{pb:euler-lagr}
\int_\Omega \frac{\nabla u \cdot \nabla v}{f(u)}=0 \qquad \forall v\in V^0=H^1_0(\Omega). 
\end{equation} 


Let $\{\Tau_h\}_{h}$ be a sequence of decompositions (meshes) of $\Omega$ into non-overlapping polygons $E$.
Each mesh $\Tau_h$ is labeled by the mesh size parameter $h$, which will be defined below, and satisfies {suitable}  regularity assumptions that are customarily made to prove the convergence of the method and derive an estimate of the approximation error. These regularity assumptions are introduced and discussed in Section~\ref{Sec:2}. 
Let $\E_h$ be the set of edges of $\Tau_h$ such that $\E_h=\E_h^i \cup \E_h^{\Gamma}$, where $\E_h^i$ and $\E_h^{\Gamma}$ are the set of interior and boundary edges, respectively. Similarly, we denote by $\V_h=\V_h^i \cup \V_h^{\Gamma}$ the set of vertices in $\Tau_h$, where $\V_h^i$ and $\V_h^{\Gamma}$ are the sets of interior and boundary vertices, respectively. Accordingly, $\V_h^E$ is the set of vertices of $E$. Moreover, $|E|$ and $|e|$ denote the area of cell $E$ and the length of edge $e$, {respectively}, $\partial E$ is the boundary of $E$, $h_E$ is the diameter of $E$ and the mesh size parameter is defined as $h=\max_{E\in\Tau_h}h_E$.  

Let us introduce the usual local lowest order conforming {V}irtual {E}lement space on the polygon $E$ (see, e.g., \cite{VEMvolley})
$$
V^E_h=\{v_h\in H^1(E): ~ \Delta v_h =0 \text{~in~} E, ~
v_h\in C^0(\partial E), v_h\vert_e \in \mathbb{P}^1(e) ~ 
\forall e \in \partial E\},
$$
{where, for $D$ $d$-dimensional domain,  $ \mathbb{P}^1(D) $ denotes the space of $d$-variate polynomials of order less than or equal to one on $D$.}
Accordingly, the global {V}irtual {E}lement space is defined as follows
$$ 
V_h^\varphi=\{v_h\in H^1(\Omega):~v_h\vert_E \in V_h^E, ~ v_h(V)=\varphi(V)~\text{~for~each~vertex~V}\in V_h^\Gamma\}.
$$
Consistently, we denote by $V_h^0$ the global VEM space with homogeneous Dirichlet boundary conditions.

Let $S^E(\cdot ,\cdot)$ be the usual stabilization term employed for constructing the VEM discretization of the Laplace problem, i.e. the Euclidean scalar product associated with the degrees of freedom (here the vertex values). See, e.g., \cite{VEMvolley, VEM:stability} for further details. Moreover, let $\Pi ^{\nabla}_E\colon V^E_{h} \rightarrow \mathbb{P}^1(E)$ the usual elliptic projection operator (see, e.g., \cite{VEMvolley}).

We introduce the local discrete function $f^E_h:V^E_h\to\mathbb{R}$ defined as
\begin{equation}
f^E_h(v_h)=\sqrt{1 + \vert \nabla \Pi^\nabla_E v_h \vert^2 + \vert E\vert ^{-1} S^E((I-\Pi^\nabla_E)v_h,(I-\Pi^\nabla_E)v_h)}.
\end{equation}
Roughly speaking, $f_h^E(\cdot)$ represents an approximation to  $\sqrt{1 + \vert (\nabla \cdot)_{\vert E} \vert^2}$. 

Having in mind the above definitions, the discrete virtual counterpart of the continuous 
minimization problem \eqref{min} reads as follows
\begin{equation}\label{VEM:minpb}
u_h = \argmin_{v_h \in V_h^\varphi} J_h(v_h),\quad \text{ with }\quad J_h(v_h) = \sum_{E\in\Tau_h} \int_E f^E_h(v_h)\,dx.
\end{equation} 
 Thus, the {V}irtual {E}lement discretization of \eqref{pb:euler-lagr} is as follows: find $u_h\in V_h^\varphi$ such that 
\begin{equation}\label{vem:pb}
A_h(u_h;u_h,v_h)=0
\end{equation}
for all $v_h\in V_h^0$, where $A_h(w_h;u_h,v_h)=\sum_{E}A_h^E(w_h;u_h,v_h)$ and 
\begin{equation}
A_h^E(w_h;u_h,v_h)=\int_E \frac{\nabla \Pi^\nabla_E u_h  \cdot  \nabla\Pi^\nabla_E v_h}{f^E_h(w_h)}\,dx + \frac{S^E((I-\Pi^\nabla_E)u_h, (I-\Pi^\nabla_E) v_h)}{f^E_h(w_h)}.
\end{equation}
 Note that as $f_h^E(w_h)$ is constant on each polygon $E$, the form $A^E_h(\cdot;\cdot,\cdot)$ can be equivalently written as 
\begin{equation}
A_h^E(w_h;u_h,v_h)=\frac{a_h^E(u_h,v_h)}{f_h^E(w_h)}
\end{equation}
where 
$$a_h^E(u_h,v_h)=\int_E 
\nabla \Pi^\nabla_E u_h \cdot  \nabla\Pi^\nabla_E v_h dx + S^E((I-\Pi^\nabla_E)u_h, (I-\Pi^\nabla_E) v_h)$$ is the classical local discrete VEM bilinear form for the Laplace problem. It is worth remembering (see, e.g., \cite{VEMvolley}) that $a_h^E(\cdot,\cdot)$ satisfies the following two crucial properties: 
  \begin{description}
  \item[$(i)$] \textbf{Consistency}: for every polynomial
    $q\in\mathbb{P}^{1}(E)$ and function $v_h\in V_h^E$ we have:
    \begin{align}
      a^{E}_h(v_h,q) = a^{E}(v_h,q);
      \label{eq:poly:r-consistency}
    \end{align}
  \item[$(ii)$] \textbf{Stability}: there exist two positive constants
    $\alpha_*$, $\alpha^*$ independent of $h$ and $E$ such that for
    every $v_h\in V_h^E$ it holds:
    \begin{align}
      \alpha_*a^{E}(v_h,v_h)\leq a^{E}_h(v_h,v_h)\leq \alpha^*a^{E}(v_h,v_h).
      \label{eq:poly:stability}
    \end{align}
  \end{description}
{ Remark that requiring that the stability condition $(ii)$ holds is equivalent to requiring that there exists positive constants $\widetilde\alpha_*$ and $\widetilde\alpha^*$ such that, for all $v_h \in V_h^E$ with $\Pi^\nabla_E v_h =0$  it holds:
	 \begin{align}
	 	{\widetilde\alpha}_*a^E(v_h,v_h)\leq
S^{E}(v_h,v_h)	 \leq {\widetilde\alpha}^*a^{E}(v_h,v_h),
	 	\label{eq:StabSE}
	 \end{align}
	 (see \cite{VEMvolley} for more details).
	}
Existence and uniqueness of the solution $u_h\in V_h^\varphi$ follow by working on the discrete cost functional $J_h (v_h)$  as in \cite{Ciarlet:book}.

For future use, {we} set $a^E(u_h,v_h)=\int_E \nabla u_h \cdot  \nabla v_h dx$.

\section{Error analysis}\label{Sec:2}

  We make the following regularity 
  assumptions on the mesh sequence $\{\Tau_h\}_h$:
\begin{itemize}
\item [(\bf{H})] there exists a constant $\rho_0>0$ independent of $\Tau_h$, such that for every element $E$ it holds:
\smallskip
\begin{enumerate}
\item[(H1)] $E$ is star-shaped with respect to all the points of a ball of radius $\rho_0 h_E$
\smallskip
\item[(H2)] every edge $e\in \E_h$ has length $\vert e \vert \geq \rho_0 h_E$.
\end{enumerate}
\end{itemize}

The assumptions (H1)-(H2) are standard (see, e.g., ~\cite{VEMvolley}) and allow to define, for every
smooth enough function $v$, an ``interpolant''  $v_I$ in 
$ V^\varphi_h$  such that it holds $|v-v_I|_{{1, \Omega}} {\lesssim} h$ (see \cite{VEMvolley}).

\newcommand{\Cuh}{C(u_h)}

We now state the main result of the paper.

{

\begin{theorem}\label{th:main}
	Let $u\in H^2(\Omega)\cap W^{1,\infty}(\Omega)$ be the continuous solution to \eqref{min}, and let
	$u_h\in V_h^\varphi$ be the VEM solution to \eqref{vem:pb}.	Letting 
	\[
	\Cuh = h^{-1} \sqrt{ \sum_E S^E((I-\Pi_E^\nabla)u_h,(I-\Pi_E^\nabla)u_h), } \]
 it holds 
	\begin{equation} | u-u_h |_{{1, \Omega}}\lesssim (1+\Cuh)^2 h.\label{apriori} \end{equation}
\end{theorem}

\begin{corollary}\label{coro:Cuh-lesssim-1}
	Assume that $\Cuh \lesssim 1$. Then it holds that
		$$ | u-u_h |_{{1, \Omega}}\lesssim  h. $$
	\end{corollary}
	
}

\begin{proof} 
By triangle inequality we have 
$$ | u-u_h|_{{1, \Omega}}\leq  | u  -u _I|_{{1, \Omega}}  
+|  u_I -u_h |_{{1, \Omega}}. $$

In the following, we adapt the ideas of \cite{Johnson-Thomee:1975} to the present context. We preliminary  observe that the following holds true
\begin{eqnarray}
| u_I- u_h|_{{1, \Omega}}&=&\left( \sum_{E\in\Tau_h} \int_E 
\frac{\vert \nabla( u_I-  u_h)\vert^2}
{f_h^E(u_h)} f^E_h(u_h)
\right)^{1/2}\nonumber \\
&\leq& 
\left( \max_E \vert f^E_h(u_h)\vert\right)^{1/2}
\left(\sum_{E\in\Tau_h} \int_E \frac{\vert \nabla( u_I- u_h)\vert^2}
{f^E_h(u_h)}
\right)^{1/2}.\label{basic}
\end{eqnarray}
The remaining part of the proof is devoted to show:
\begin{itemize}
\item[\emph{(i)}] $\displaystyle{\left(\sum_{E\in\Tau_h} \int_E \frac{\vert \nabla( u_I- u_h)\vert^2}
{f_h^E(u_h)}
\right)^{1/2}\lesssim (1+ \Cuh) h}$;
\item[\emph{(ii)}] $ \displaystyle{\max_E \vert f^E_h(u_h) {\vert} \lesssim (1+\Cuh)^2}$.
\end{itemize}

Let us first prove $(i)$. We start by observing that, thanks to \eqref{eq:StabSE}, we have
\begin{equation}\label{aux:silvia}
| u_h - \Pi^\nabla_E u_h |_{1,h} \lesssim \Cuh h
\end{equation}
where {where $| v|_{1,h}^2=\sum_{E\in\Tau_h}{\| \nabla v\|_{0,E}^2}$.}
By using the stability property of $a_h^E(\cdot,\cdot)$, as $f_h^E(u_h)$ is constant on $E$, we get the following inequalities with $\delta_h= u_h- u_I$
\begin{eqnarray}
\sum_{E\in\Tau_h} \int_E \frac{\vert \nabla( u_I- u_h)\vert^2}
{f_h^E(u_h)} &\lesssim & 
\sum_{E\in\Tau_h} \frac{a_h^E( \delta_h, \delta_h)}
{f_h^E(u_h)} \nonumber\\
&\lesssim& \left\vert -\sum_{E\in\Tau_h} \frac{a_h^E(  u_I, \delta_h)}{f_h^E(u_h)}\right\vert,
\label{aux:1}
\end{eqnarray}
where in the last step we employ \eqref{vem:pb} with $\delta_h\in V_h^0$. Let $ u_\pi\vert_E$ be the $L^2(E)$ projection of $u$ onto $\mathbb{P}^1(E)$. By employing the consistency and stability properties of $a_h^E (\cdot,\cdot)$ together with the fact that $u$ is solution to \eqref{min}, it is immediate to check that the following holds

\begin{eqnarray}
-\sum_{E\in\Tau_h} \frac{a_h^E(  u_I,\delta_h)}{ f^E_h(u_h)}
&=&
-\sum_{E\in\Tau_h} \Bigg\{
\frac{a_h^E( u_I- u_\pi,\delta_h)}{ f^E_h(u_h)}
+
\frac{a_h^E(  u_\pi, \delta_h)}{ f^E_h(u_h)} \\ &&
\pm 
\frac{a^E(  u, \delta_h)}{ f^E_h(u_h)}
-
\int_E\frac{\nabla u  \cdot  \nabla\delta_h}{f(u)}dx
\Bigg\}\nonumber\\
&=&
\sum_{E\in\Tau_h} 
\frac{a_h^E(  u_\pi- u_I,\delta_h)}{ f^E_h(u_h)}
+
\sum_{E\in\Tau_h} 
\frac{a^E( u- u_\pi,\delta_h)}{ f^E_h(u_h)}\nonumber\\
&&+
\sum_{E\in\Tau_h} 
\int_E \nabla u \cdot  \nabla \delta_h\left(\frac{1}{f(u)}-\frac{1}{ f^E_h(u_h)}\right)dx\nonumber\\[4mm]
&=&A+B+C. 
\end{eqnarray}
We now bound the three terms separately. 
By combining the Cauchy-Schwarz inequality with the fact that $f_h^E(u_h)$ is constant and larger than $1$ on each polygon $E$, we have 
\begin{eqnarray}
A\leq (| u- u_\pi|_{1,h}+| u- u_I|_{{1, \Omega}})\left(\sum_{E\in\Tau_h} \int_E \frac{\vert \nabla( u_I- u_h)\vert^2}
{ f^E_h(u_h)}dx
\right)^{1/2},
\end{eqnarray}
and
\begin{eqnarray}
B\leq | u- u_\pi|_{1,h}\left(\sum_{E\in\Tau_h} \int_E \frac{\vert \nabla( u_I- u_h)\vert^2}
{ f^E_h(u_h)}dx
\right)^{1/2}.
\end{eqnarray}
Finally, setting 
$\gamma=\max_{\bar{\Omega}} \frac{\vert \nabla u\vert}{f(u)}$, employing the definitions of $f(\cdot)$ and $f^
E_h(\cdot)$ and observing that $f(u)\geq |\nabla u|$,  the following holds
\begin{eqnarray}
C&=& \sum_{E\in \Tau_h} \int_E \nabla u \cdot  \nabla \delta_h
\frac{[f_h^E(u_h)]^2-f^2(u)}{f(u)f_h^E(u_h)(f(u)+f_h^E(u_h))}dx\nonumber\\ 
&\leq&  \gamma\sum_{E\in \Tau_h} \int_E  | \nabla \delta_h |
\frac{| f^2(u)-[f_h^E(u_h)]^2 |}{f_h^E(u_h)(f(u)+f_h^E(u_h))}dx\nonumber\\ 
&=& \gamma \sum_{E\in \Tau_h}\Bigg\{\int_E | \nabla \delta_h |
\frac{ \left| |\nabla u|^2-|\nabla \Pi_E^\nabla u_h|^2 \right| }{f_h^E(u_h)(f(u)+f_h^E(u_h))}dx\nonumber\\ &&+ 
 \int_E | \nabla \delta_h | 
 \frac{
 	|E|^{-1} S^E((I-\Pi_E^\nabla)u_h,(I-\Pi_E^\nabla)u_h))}{f_h^E(u_h)(f(u)+f_h^E(u_h))}dx \Bigg\}
\nonumber =  C.I + C.II \nonumber
\end{eqnarray}
As $f_h^E(u_h)\geq |\nabla\Pi_E^\nabla u_h|$, we can bound
\begin{eqnarray*}
C.I &=& \gamma \sum_{E\in \Tau_h} \int_E | \nabla \delta_h |
\frac{\left| |\nabla u|^2-|\nabla \Pi_E^\nabla u_h|^2 \right| }{f_h^E(u_h)(f(u)+f_h^E(u_h))}dx\\ &\leq& \gamma \sum_{E\in \Tau_h} \int_E | \nabla \delta_h |
\frac{|\nabla (u -\Pi_E^\nabla u_h)|(|\nabla u|+|\nabla \Pi_E^\nabla u_h|)
	}{f_h^E(u_h)(f(u)+f_h^E(u_h))}dx \\&\leq& \gamma \sum_{E\in \Tau_h} \int_E | \nabla \delta_h | 
	\frac{|\nabla (u -\Pi_E^\nabla u_h)|}{f_h^E(u_h)}dx.
\end{eqnarray*}


Now, employing the Cauchy-Schwarz inequality and noticing that $f_h^E(u_h)\geq 1$,  we have the following
\begin{eqnarray}
C.I&\lesssim&\gamma \left(\sum_{E\in\Tau_h} \int_E \frac{\vert \nabla \delta_h\vert^2}
{ f^E_h(u_h)}dx\right)^{1/2} \left\{
\left(\sum_{E\in\Tau_h} \int_E \frac{\vert \nabla\delta_h\vert^2}
{ f^E_h(u_h)}dx\right)^{1/2}
+ 
| u - u_I|_{{1, \Omega}} + | u_h- \Pi_E^\nabla u_h  |_{1,h}\right\}\nonumber \\
&\leq&
\gamma \left(\sum_{E\in\Tau_h} \int_E \frac{\vert \nabla \delta_h\vert^2}
{ f^E_h(u_h)}dx\right)^{1/2} \left\{
\left(\sum_{E\in\Tau_h} \int_E \frac{\vert \nabla\delta_h\vert^2}
{ f^E_h(u_h)}dx\right)^{1/2}
+ 
| u - u_I|_{{1, \Omega}} + \Cuh h\right\},
\end{eqnarray}
where we used the stability property \eqref{aux:silvia} and the definition of the constant $C(u_h)$.
On the other hand, as $f_h^E(u_h) > 1$ clearly implies $[f_h^E(u_h)]^2\geq [f_h^E(u_h)]^{3/2}$, we have
\begin{eqnarray}
C.II &\leq & \gamma \sum_{E\in \Tau_h} 
\int_E\frac{|\nabla\delta_h|}{2[f_h^E(u_h)]^{1/2}} \frac{|E|^{-1} S^E((I-\Pi_E^\nabla)u_h,(I-\Pi_E^\nabla)u_h))}{f_h^E(u_h)}dx\nonumber\\ &\leq&
\gamma \sum_{E\in \Tau_h} 
\int_E\frac{|\nabla\delta_h|}{2[f_h^E(u_h)]^{1/2}} |E|^{-1/2} (S^E((I-\Pi_E^\nabla)u_h,(I-\Pi_E^\nabla)u_h)))^{1/2}dx\nonumber\\ &\leq&
\gamma \left(\sum_{E\in\Tau_h} \int_E \frac{\vert \nabla\delta_h\vert^2}
{ f^E_h(u_h)}\right)^{1/2}  \Cuh h,
\end{eqnarray}
where we used $f_h^E(u_h)\geq  |E|^{-1/2} (S^E((I-\Pi_E^\nabla)u_h,(I-\Pi_E^\nabla)u_h)))^{1/2} $ and employed the Cauchy-Schwarz inequality once  again.
Setting $$T=\sum_{E\in\Tau_h} \int_E \frac{\vert \nabla \delta_h \vert^2}
{ f^E_h(u_h)}dx$$ and plugging the above inequalities for $A$, $B$, $C$ into \eqref{aux:1} we obtain 
\begin{multline*}
T\lesssim T^\frac 1 2 (| u- u_\pi|_{1,h}+| u- u_I|_{{1, \Omega}}) 
+
 T^\frac 1 2 | u- u_\pi|_{1,h}
\\+
\gamma T^\frac 1 2 \left(T^\frac 1 2 + 
| u - u_I|_{{1, \Omega}}
+
2C(u_h) h
\right). 
\end{multline*}
Noticing that $\gamma<1$ we get 
\begin{equation}
T^{1/2}\lesssim \frac{1}{1-\gamma}( | u- u_\pi|_{1,h} +  | u- u_I|_{1} + \Cuh h  ),
\end{equation}
which, using standard error estimates, implies $T^{1/2} \lesssim (1+C(u_h)) h $.

\medskip
Finally, we  prove (ii).  In particular, from (i) we have 
$$
\left(\int_E \frac{\vert \nabla\delta_h\vert^2}
{ f^E_h(u_h)}dx
\right)^{1/2}\lesssim(1+\Cuh) h
$$ 
{for any $E\in \Tau_h$,}
which implies
\begin{eqnarray}
\left(\int_E \frac{\vert \nabla u_h\vert^2}
{ f^E_h(u_h)}dx
\right)^{1/2}&\lesssim& 
\left(\int_E \frac{\vert \nabla \delta_h\vert^2}
{ f^E_h(u_h)}dx
\right)^{1/2}
+ 
\left(\int_E \frac{\vert \nabla u_I \vert^2}
{ f^E_h(u_h)}dx
\right)^{1/2}\nonumber\\
&\lesssim&(1+\Cuh) h + |u|_{W^{1,\infty}} \left (\int_E dx\right)^{1/2} \nonumber\\
&\lesssim& (1+\Cuh) h, \label{aux:2}
\end{eqnarray}
where we employed the fact that $f_h^E(u_h)\geq 1$ on each $E$, the $H^1$-stability of the  interpolation operator $(\cdot)_I$ and $|E|\simeq h^2$.

On the other hand, using the fact that $f_h^E(u_h)$ is constant on each $E$ and employing the $H^1$-orthogonality property of the elliptic projector $\Pi^\nabla_E$  we have 
\begin{eqnarray}
\int_E \frac{\vert \nabla u_h\vert^2}
{ f^E_h(u_h)} dx &=& \int_E \frac{\vert \nabla \Pi^\nabla_E u_h \vert^2}
{ f^E_h(u_h)} dx + \int_E \frac{\vert \nabla (I-\Pi_E^\nabla)u_h\vert^2}
{ f^E_h(u_h)} dx\nonumber\\
&\gtrsim &
\int_E \frac{\vert \nabla \Pi^\nabla_E u_h \vert^2}
{ f^E_h(u_h)} dx
+
\frac{S^E((I-\Pi_E^\nabla)u_h,(I-\Pi_E^\nabla)u_h)}{ f^E_h(u_h)}{,}
\label{aux:3}
\end{eqnarray}
where in the last step we employed \eqref{eq:StabSE}. 
Combining \eqref{aux:2} and \eqref{aux:3}, and observing that $\Pi^\nabla_E u_h$ and $S^E((I-\Pi_E^\nabla)u_h,(I-\Pi_E^\nabla)u_h)$ are both constant on $E$ yield
$$\int_E \frac{\vert \nabla \Pi^\nabla_E u_h\vert^2+
\vert E \vert ^{-1} S^E((I-\Pi_E^\nabla)u_h,(I-\Pi_E^\nabla)u_h)}
{ f^E_h(u_h)}dx\lesssim(1+\Cuh)^2 h^2{,}$$
and thus
$$\frac{\vert \nabla \Pi^\nabla_E u_h\vert^2+
\vert E \vert ^{-1} S^E((I-\Pi_E^\nabla)u_h,(I-\Pi_E^\nabla)u_h)}
{ f^E_h(u_h)}\lesssim (1+\Cuh)^2{,}$$
which, recalling the definition of $f_h^E(u_h)$,  implies $$\vert \nabla \Pi^\nabla_E u_h\vert^2+
\vert E \vert ^{-1} S^E((I-\Pi_E^\nabla)u_h,(I-\Pi_E^\nabla)u_h)\lesssim (1+\Cuh)^4.$$ This yields (ii).
By combining (i) and (ii) with \eqref{basic} we finally obtain the thesis.
\end{proof}

{
	\begin{remark}
		Observe that, while \eqref{apriori} is not properly an {\em a priori} estimate on the error, as the quantity $\Cuh$ on the right hand side depends on the discrete solution and, consequently on $h$, such a quantity can be computed {\em a posteriori}, allowing us to check whether it remains bounded, thus providing a useful bound. Observe also that such a quantity is obtained by combining local contributions, so that, should it be too big, its distribution might (heuristically) provide some information on how to refine the mesh in order to obtain a better solution. 
	\end{remark}

	}

\section{Numerical Experiments}\label{Sec:3}

The discrete VE problem \eqref{vem:pb} is solved using a classical fixed point algorithm, i.e. iterate on $k$ the following: given $u_h^k \in V_h^\varphi$, find $u_h^{k+1} \in V_h^\varphi$ such that 
\[
A_h(u_h^k;u_h^{k+1}, v_h) = 0\qquad\forall\,v_h\in V_h^0  \quad (\textit{linearized~problem}).
\]
Fixed point iterations are stopped as soon as $||u_h^{k+1} - u_h^k||_\infty/||u_h^k||_\infty$ is less than a prescribed tolerance  $\mathtt{tol} = \num[retain-unity-mantissa = false]{1e-9}$, whereas at each iteration, the discrete linear system is solved using 
a direct solver.

\begin{figure}
\centering
\subfloat{\includegraphics[width=0.25\textwidth]{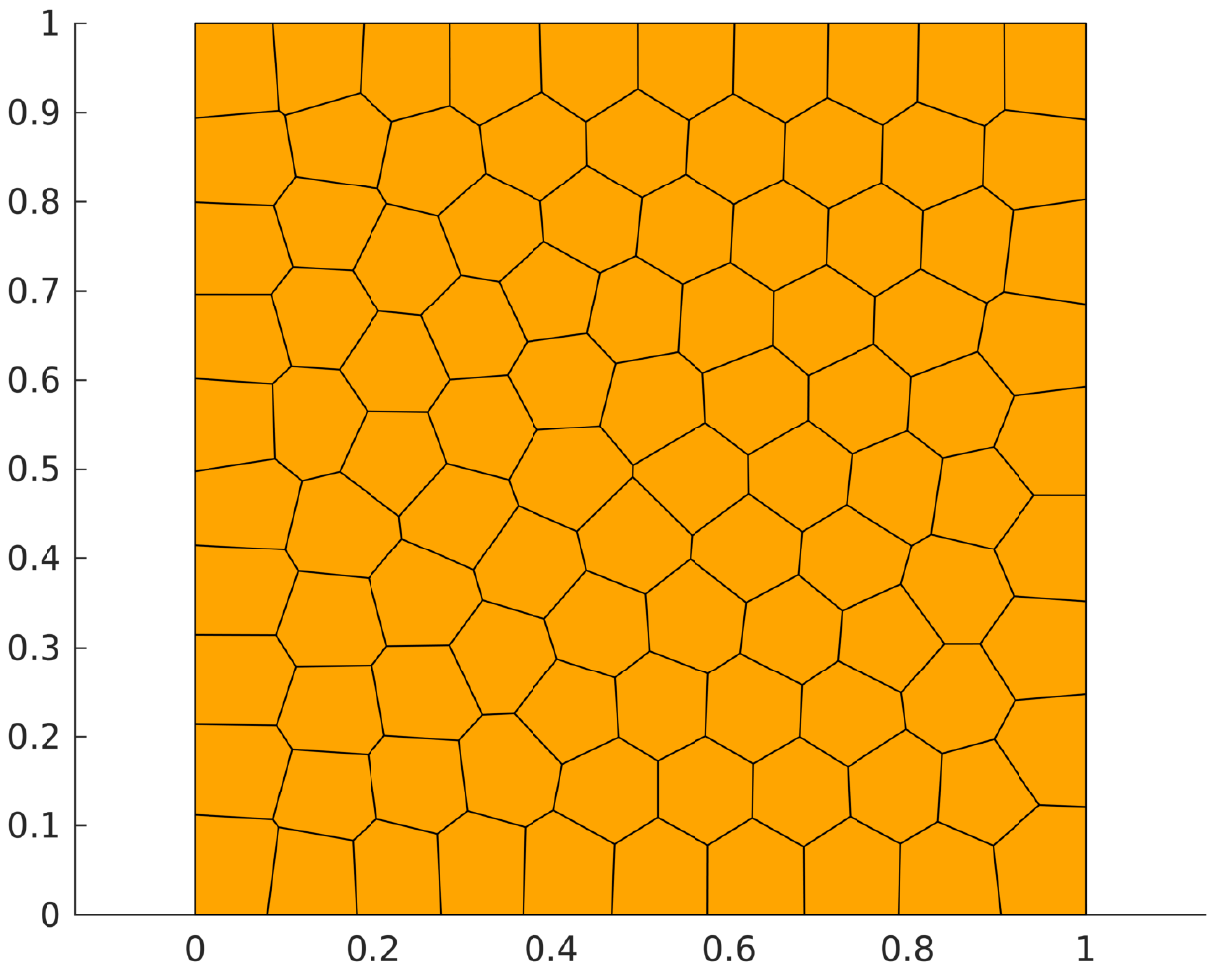}}
\subfloat{\includegraphics[width=0.25\textwidth]{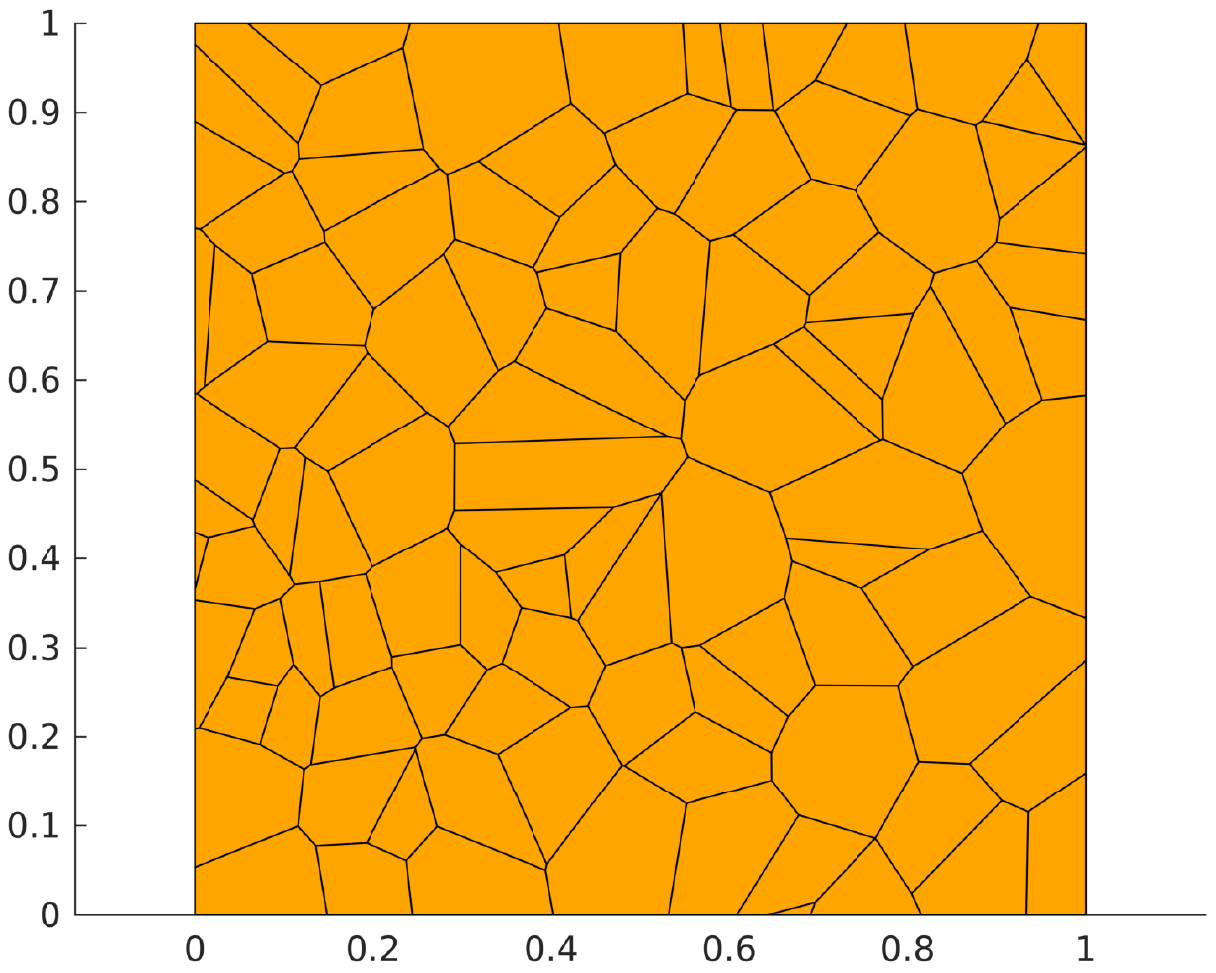}}
\caption{{Example of the meshes used in the numerical tests.}}
\label{fig:example-meshes}
\end{figure}

To assess the convergence properties of our Virtual Element discretization, we introduce the following error quantities:
\begin{align*}
e_{H^1} &= \frac{||\nabla u - \Pi^0_0\nabla u_h||_{L^2(\Omega)}}{||\nabla u||_{L^2(\Omega)}}, & e_{L^2} &= \frac{||u - \Pi^0_1 u_h||_{L^2(\Omega)}}{||u||_{L^2(\Omega)}},
\end{align*}
where $\Pi^0_k$ is the $L^2$-projection onto the space of polynomials of degree $k$, {$k=0,1$}. The exact solution $u$ is evaluated analytically, whenever possible. Otherwise, it is approximated by the solution $u_h^\text{FEM}$ computed with the finite element method on a {very} fine grid of $\Omega$. Estimated convergence rates (ecr) are computed with respect to the total number of degrees of freedom $N$, under the assumption $N \approx O(h^{-2})$.
All {the} numerical experiments are performed on Voronoi meshes that are either uniform or random, see Figure~\ref{fig:example-meshes}.
For each mesh, we collect  the following informations (see tables below):
\begin{itemize}
\item the maximum diameter over all the elements of the mesh {($h$)};
\item the number of degrees of freedom {($N$)};
\item the number of fixed-point iterations required to reach convergence {(It)};
\item {the computed errors} $e_{H^1}$ and $e_{L^2}$ {measured in the $H^1$ and $L^2$ norms, respectively,} and the corresponding estimated convergence rates {(ecr)};
\item the constant $\Cuh$ defined in Theorem~\ref{th:main}, computed by using either the actual diameter $h$ (see the column named $C_1$) or $1/\sqrt{N} \approx h$ (see the column named $C_2$).
\end{itemize}

\subsection{Test 1}\label{sec:concus}
\begin{figure}
\centering
\subfloat{\includegraphics[width=0.5\textwidth]{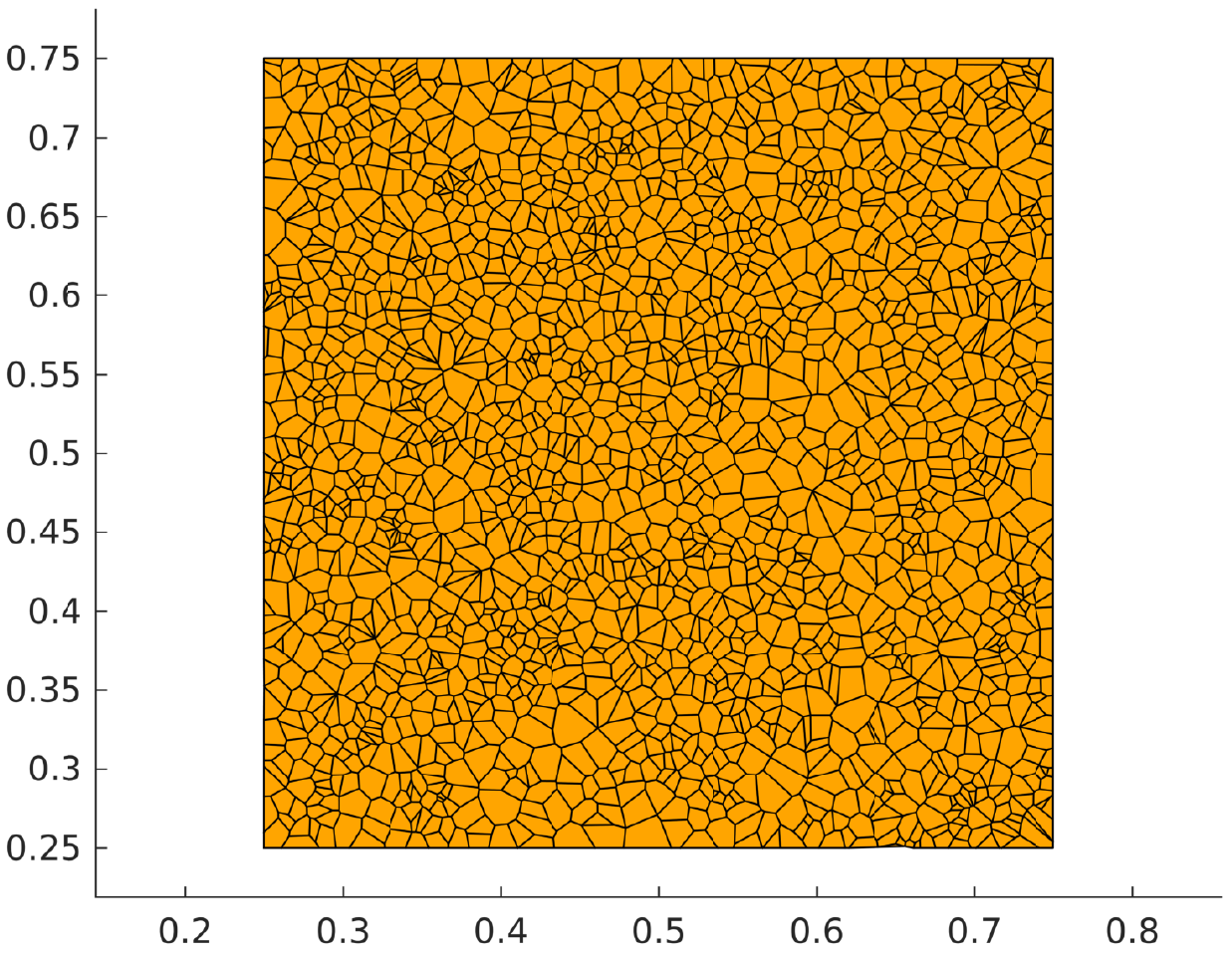}}
\subfloat{\includegraphics[width=0.5\textwidth]{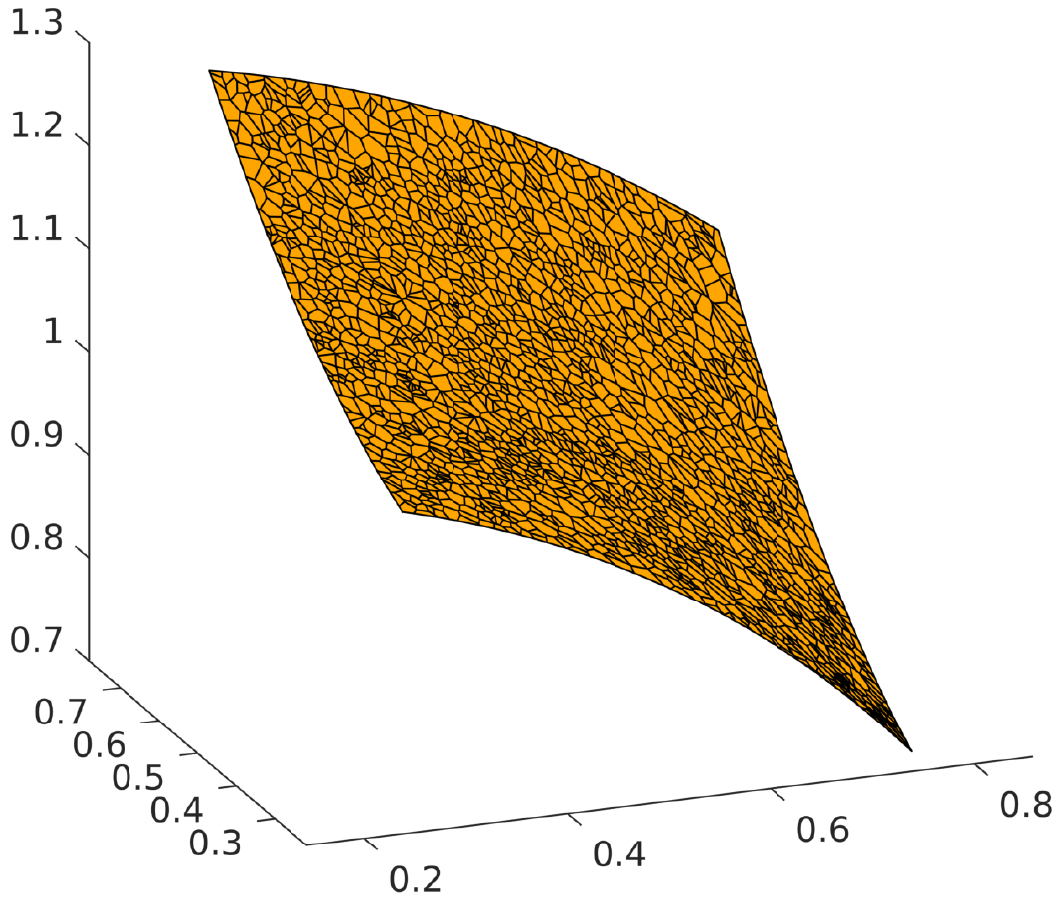}}
\caption{{Example~\ref{sec:concus}: example of the computational mesh (left) and corresponding computed solution (right).}}
\label{fig:example-concus}
\end{figure}

Here we consider a test problem originally proposed by Concus~\cite{consus1966} 
that provided the following analytic solution to the minimal surface problem on the square $\Omega = (0.25, 0.75)\times(0.25, 0.75)$:
\[
u(x,y) = \sqrt{\cosh^2(y) - x^2}.
\]
Note that $u \in H^2(\Omega)\cap W^{1,\infty}(\Omega)$. An example of computed solution  on a coarse mesh is shown in Figure~\ref{fig:example-concus}. Experiments are performed on uniform (Table~\ref{tab:concus-uniformVoro}) and random Voronoi meshes (Table~\ref{tab:concus-randVoro}). The assumption $\Cuh \lesssim 1$ is verified and the rate of convergence in the $H^1$-norm is in agreement with Theorem  \ref{th:main}. Moreover, the reported rate of convergence in the $L^2$-norm seems to be $2$.
\begin{table}
\centering
\caption{{Example~\ref{sec:concus}: computed errors and estimated convergence rates (uniform Voronoi meshes).}}
\label{tab:concus-uniformVoro}
\begin{tabular}{
c
S[table-format=1.{\roundPrecision}e-1]
S[table-format=6.0]
S[table-format=2.0]
S[table-format=1.{\roundPrecision}e-1]
S[table-format=1.{\roundPrecision}]
S[table-format=1.{\roundPrecision}e-1]
S[table-format=1.{\roundPrecision}]
S[table-format=1.{\roundPrecision}]
S[table-format=1.{\roundPrecision}]
}
\toprule
{Mesh} & {$h$} & {$N$} & {It} & {$e_{H^1}$} & {ecr} & {$e_{L^2}$} & {ecr} & {$C_1$} & {$C_2$}\\
\midrule
u-concus$_{1}$   &   1.820473e-02   &   4082   &   17   &   7.315457e-03   &   {-}   &   1.228661e-05   &   {-}   &   0.1337312   &   0.1555441\\
u-concus$_{2}$   &   1.233497e-02   &   8165   &   17   &   5.138253e-03   &   1.019160   &   6.175387e-06   &   1.984619   &   0.1389607   &   0.1548845\\
u-concus$_{3}$   &   8.472495e-03   &   16323   &   17   &   3.640650e-03   &   0.994779   &   3.206931e-06   &   1.891843   &   0.1457149   &   0.1577304\\
u-concus$_{4}$   &   6.358370e-03   &   32657   &   17   &   2.564317e-03   &   1.010751   &   1.521949e-06   &   2.149502   &   0.1365472   &   0.1568977\\
u-concus$_{5}$   &   4.512142e-03   &   65293   &   17   &   1.817913e-03   &   0.993043   &   7.959064e-07   &   1.871367   &   0.1371768   &   0.1581600\\
u-concus$_{6}$   &   3.002593e-03   &   130567   &   17   &   1.283903e-03   &   1.003706   &   4.037045e-07   &   1.959010   &   0.1457675   &   0.1581517\\
u-concus$_{7}$   &   2.235510e-03   &   261206   &   17   &   9.077853e-04   &   0.999828   &   1.979326e-07   &   2.055762   &   0.1384573   &   0.1581918\\
u-concus$_{8}$   &   1.522826e-03   &   522279   &   17   &   6.420748e-04   &   0.999587   &   9.936792e-08   &   1.989045   &   0.1440562   &   0.1585381\\
\bottomrule
\end{tabular}
\end{table}
\begin{table}
\centering
\caption{{Example~\ref{sec:concus}: computed errors and estimated convergence rates (random Voronoi meshes).}}
\label{tab:concus-randVoro}
\begin{tabular}{
c
S[table-format=1.{\roundPrecision}e-1]
S[table-format=6.0]
S[table-format=2.0]
S[table-format=1.{\roundPrecision}e-1]
S[table-format=1.{\roundPrecision}]
S[table-format=1.{\roundPrecision}e-1]
S[table-format=1.{\roundPrecision}]
S[table-format=1.{\roundPrecision}]
S[table-format=1.{\roundPrecision}]
}
\toprule
{Mesh} & {$h$} & {$N$} & {It} & {$e_{H^1}$} & {ecr} & {$e_{L^2}$} & {ecr} & {$C_1$} & {$C_2$}\\
\midrule
concus$_{1}$   &   3.406596e-02   &   3716   &   17   &   9.316064e-03   &   {-}   &   2.736716e-05   &   {-}   &   0.09103010   &   0.1890356\\
concus$_{2}$   &   2.526368e-02   &   7450   &   17   &   6.407736e-03   &   1.076057   &   1.114473e-05   &   2.583154   &   0.08371536   &   0.1825491\\
concus$_{3}$   &   1.717608e-02   &   14693   &   17   &   4.582240e-03   &   0.987453   &   5.944094e-06   &   1.851026   &   0.08905751   &   0.1854171\\
concus$_{4}$   &   1.293438e-02   &   29487   &   17   &   3.223429e-03   &   1.009915   &   3.044994e-06   &   1.920529   &   0.08353030   &   0.1855261\\
concus$_{5}$   &   9.161274e-03   &   59011   &   17   &   2.263267e-03   &   1.019456   &   1.581828e-06   &   1.887984   &   0.08304098   &   0.1848055\\
concus$_{6}$   &   7.240372e-03   &   118053   &   17   &   1.601315e-03   &   0.997921   &   7.453749e-07   &   2.170287   &   0.07448809   &   0.1853046\\
concus$_{7}$   &   4.918414e-03   &   235898   &   17   &   1.134128e-03   &   0.996615   &   3.709904e-07   &   2.015732   &   0.07754201   &   0.1852355\\
concus$_{8}$   &   3.568262e-03   &   472263   &   17   &   8.009051e-04   &   1.002330   &   1.840569e-07   &   2.019577   &   0.07556860   &   0.1853062\\
\bottomrule
\end{tabular}
\end{table}


\subsection{Test 2}\label{sec:catenoid}
\begin{figure}
\centering
\subfloat{\includegraphics[width=0.5\textwidth]{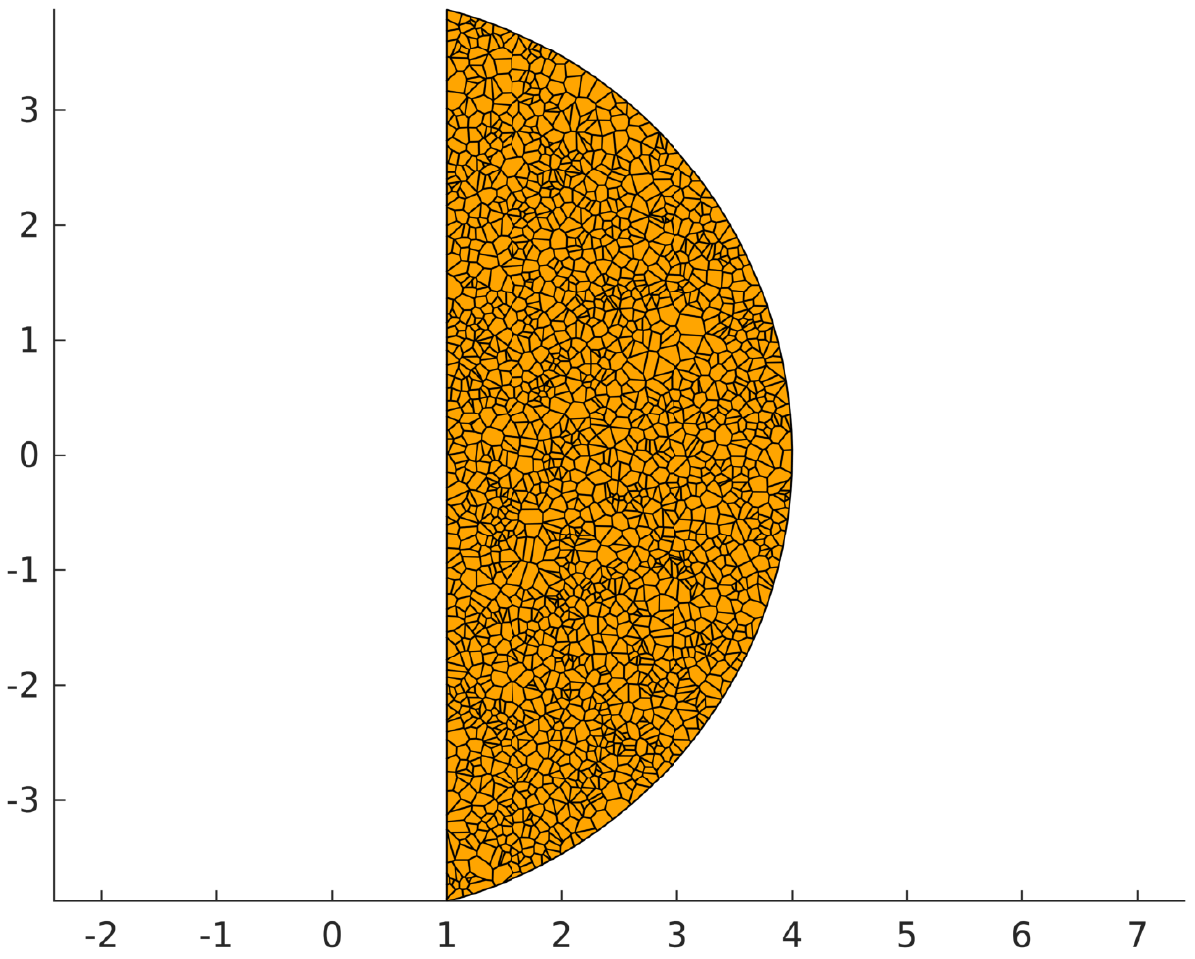}}
\subfloat{\includegraphics[width=0.5\textwidth]{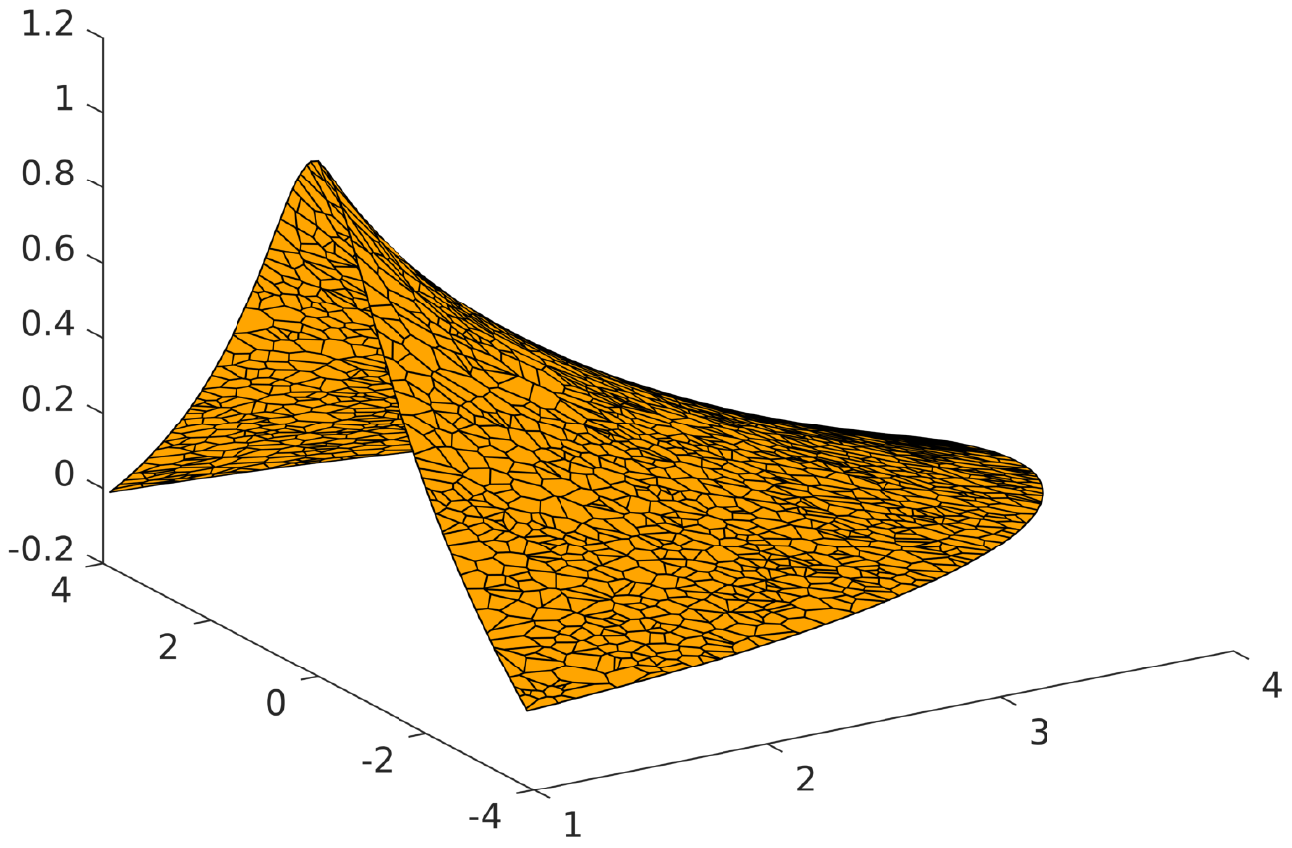}}
\caption{{Example~\ref{sec:catenoid}: example of the computational mesh (left) and corresponding computed solution (right).}}
\label{fig:example-catenoid}
\end{figure}

Here we consider another test problem for which an analytic solution is known~\cite{nitsche1965}. Let us consider the following convex domain
\[
\Omega = \Set{\mathbf x = (x, y)\in\mathbb R^2 | \lVert\mathbf x\rVert_2 < 4 \cap x > 1}.
\]
An explicit example of minimal surface on $\Omega$ is given by
\[
u(x,y) = a\log\left(\frac{b + \sqrt{b^2 - a^2}}{r + \sqrt{r^2 - a^2}}\right),
\]
where we take $a = 0.75, b = 4$ and $r = \sqrt{x^2 + y^2}$. Note that $u\in H^2(\Omega)\cap W^{1,\infty}(\Omega)$. This minimal surface is also known as \emph{catenoid}. A typical solution on a coarse mesh is shown in Figure~\ref{fig:example-catenoid}. Experiments are performed on uniform Voronoi meshes (Table~\ref{tab:catenoid-uniformVoro}) and random Voronoi meshes (Table~\ref{tab:catenoid-randVoro}). Again, $\Cuh \lesssim 1$ and the rate of convergence in the $H^1$-norm is in agreement with Theorem \ref{th:main}, whereas the computed rate of convergence in the $L^2$-norm seems to be $2$.
\begin{table}
\centering
\caption{{Example~\ref{sec:catenoid}: computed errors and estimated convergence rates (uniform Voronoi meshes).}}
\label{tab:catenoid-uniformVoro}
\begin{tabular}{
c
S[table-format=1.{\roundPrecision}e-1]
S[table-format=6.0]
S[table-format=2.0]
S[table-format=1.{\roundPrecision}e-1]
S[table-format=1.{\roundPrecision}]
S[table-format=1.{\roundPrecision}e-1]
S[table-format=1.{\roundPrecision}]
S[table-format=1.{\roundPrecision}]
S[table-format=1.{\roundPrecision}]
}
\toprule
{Mesh} & {$h$} & {$N$} & {It} & {$e_{H^1}$} & {ecr} & {$e_{L^2}$} & {ecr} & {$C_1$} & {$C_2$}\\
\midrule
u-sector$_{1}$   &   1.408847e-01   &   4080   &   20   &   2.571832e-02   &   {-}   &   3.745642e-04   &   {-}   &   0.198131   &   1.782979\\
u-sector$_{2}$   &   1.031659e-01   &   8158   &   21   &   1.748529e-02   &   1.113704   &   1.720191e-04   &   2.246083   &   0.186410   &   1.736990\\
u-sector$_{3}$   &   6.978080e-02   &   16309   &   21   &   1.262748e-02   &   0.939731   &   8.938975e-05   &   1.889945   &   0.198421   &   1.768227\\
u-sector$_{4}$   &   5.177668e-02   &   32640   &   22   &   8.940494e-03   &   0.995312   &   4.514633e-05   &   1.969086   &   0.187909   &   1.757746\\
u-sector$_{5}$   &   3.587076e-02   &   65271   &   22   &   6.378085e-03   &   0.974656   &   2.291350e-05   &   1.957211   &   0.194161   &   1.779358\\
u-sector$_{6}$   &   2.557910e-02   &   130572   &   23   &   4.431724e-03   &   1.050163   &   1.115385e-05   &   2.076624   &   0.188380   &   1.741181\\
u-sector$_{7}$   &   1.807119e-02   &   261077   &   23   &   3.140899e-03   &   0.993749   &   5.583417e-06   &   1.997382   &   0.190167   &   1.755928\\
u-sector$_{8}$   &   1.296064e-02   &   522210   &   23   &   2.216130e-03   &   1.006114   &   2.763452e-06   &   2.029040   &   0.187993   &   1.760722\\
\bottomrule
\end{tabular}
\end{table}
\begin{table}
\centering
\caption{{Example~\ref{sec:catenoid}: computed errors and estimated convergence rates (random Voronoi meshes).}}
\label{tab:catenoid-randVoro}
\begin{tabular}{
c
S[table-format=1.{\roundPrecision}e-1]
S[table-format=6.0]
S[table-format=2.0]
S[table-format=1.{\roundPrecision}e-1]
S[table-format=1.{\roundPrecision}]
S[table-format=1.{\roundPrecision}e-1]
S[table-format=1.{\roundPrecision}]
S[table-format=1.{\roundPrecision}]
S[table-format=1.{\roundPrecision}]
}
\toprule
{Mesh} & {$h$} & {$N$} & {It} & {$e_{H^1}$} & {ecr} & {$e_{L^2}$} & {ecr} & {$C_1$} & {$C_2$}\\
\midrule
sector$_{1}$   &   2.783233e-01   &   4198   &   19   &   3.037020e-02   &   {-}   &   6.536296e-04   &   {-}   &   0.121057   &   2.183039\\
sector$_{2}$   &   2.262086e-01   &   8330   &   19   &   2.303112e-02   &   0.807336   &   3.991139e-04   &   1.439738   &   0.109930   &   2.269595\\
sector$_{3}$   &   1.475084e-01   &   16588   &   21   &   1.536174e-02   &   1.175834   &   1.606314e-04   &   2.642606   &   0.114747   &   2.179986\\
sector$_{4}$   &   1.074039e-01   &   33080   &   22   &   1.106745e-02   &   0.950009   &   8.749437e-05   &   1.760343   &   0.113671   &   2.220512\\
sector$_{5}$   &   7.993085e-02   &   65973   &   22   &   7.805126e-03   &   1.011790   &   4.463371e-05   &   1.950076   &   0.108630   &   2.230221\\
sector$_{6}$   &   5.475817e-02   &   131673   &   23   &   5.477103e-03   &   1.025081   &   2.080637e-05   &   2.208816   &   0.111298   &   2.211494\\
sector$_{7}$   &   4.037713e-02   &   262975   &   23   &   3.897524e-03   &   0.983711   &   1.057072e-05   &   1.957886   &   0.108441   &   2.245369\\
sector$_{8}$   &   2.853023e-02   &   525468   &   23   &   2.724878e-03   &   1.034100   &   5.122127e-06   &   2.093286   &   0.107732   &   2.228042\\
\bottomrule
\end{tabular}
\end{table}


\subsection{Test 3}\label{sec:scherk}
\begin{figure}
\centering
\subfloat{\includegraphics[width=0.5\textwidth]{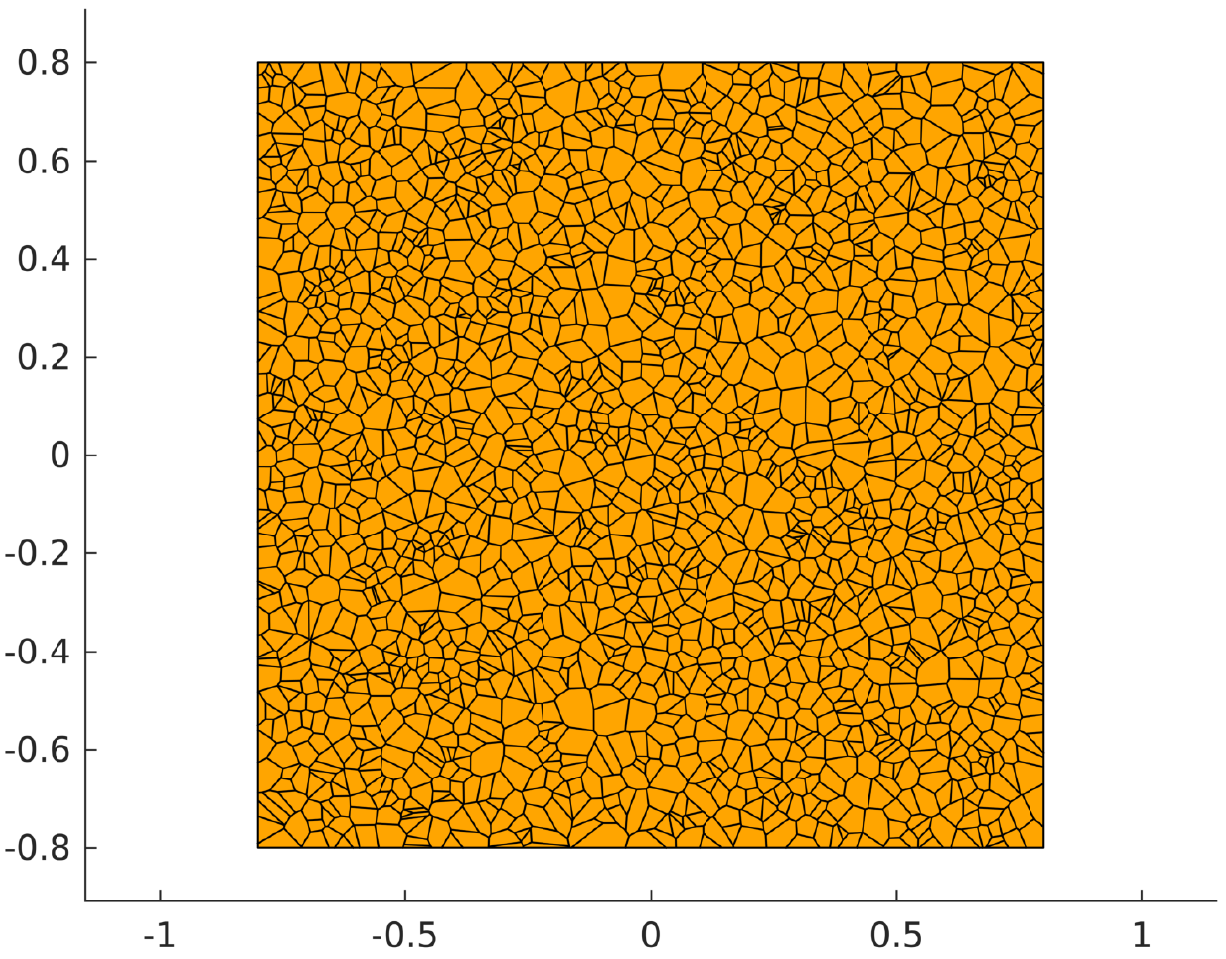}}
\subfloat{\includegraphics[width=0.5\textwidth]{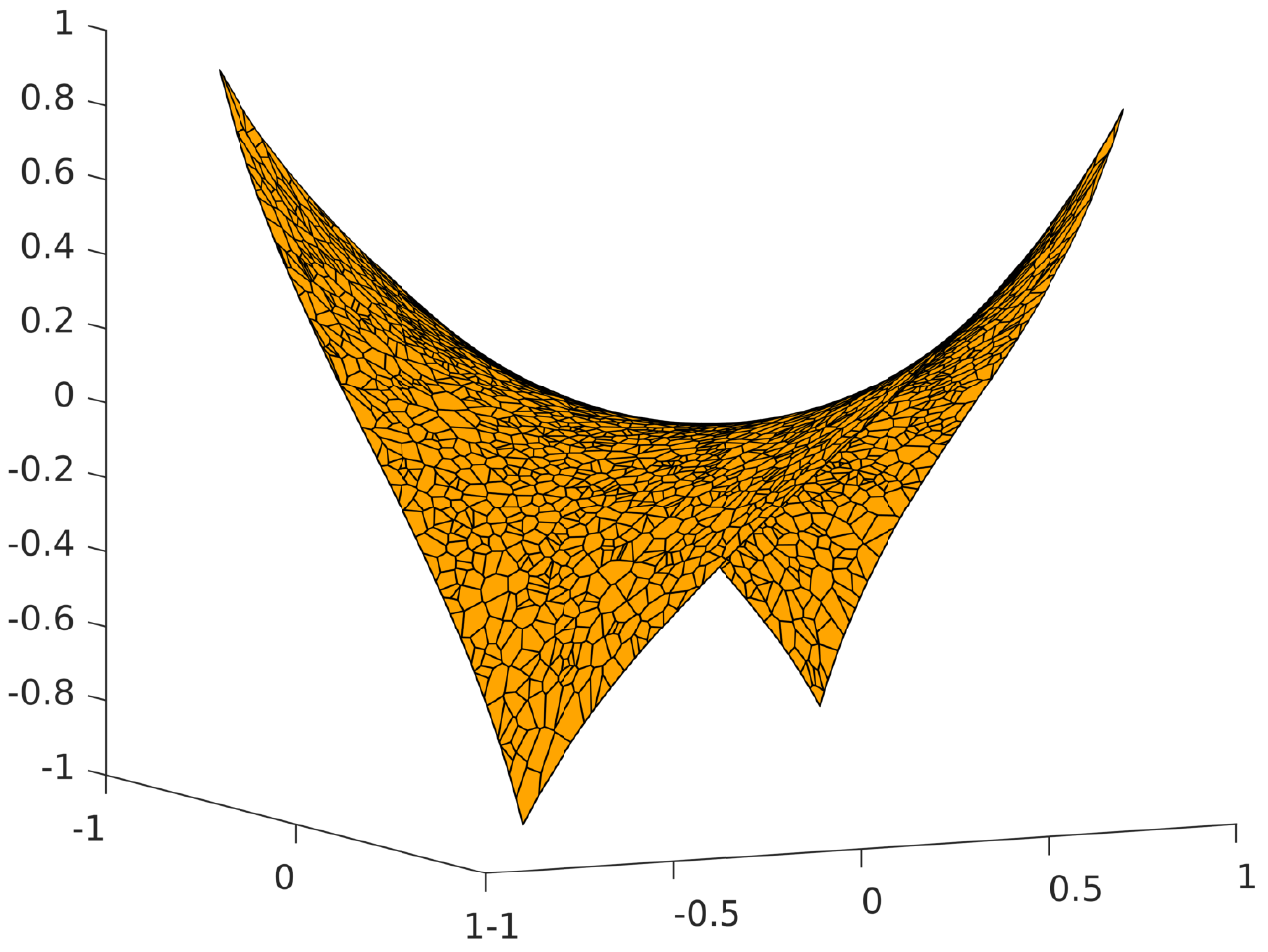}}
\caption{{Example~\ref{sec:scherk}: example of the computational mesh (left) and corresponding computed solution (right).}}
\label{fig:example-scherk}
\end{figure}

Here we consider the so called Scherk's fifth surface~\cite{trasdal2011} which is another minimal surface that can be expressed on $\Omega = (-0.8,0.8)\times(-0.8,0.8)$ as follows
\[
u(x,y) = \sin^{-1}( \sinh x \sinh y ).
\]
A typical solution on a coarse mesh is shown in Figure~\ref{fig:example-scherk}. Experiments are performed on uniform Voronoi meshes (Table~\ref{tab:scherk-uniformVoro}) and random Voronoi meshes (Table~\ref{tab:scherk-randVoro}). The assumption $\Cuh \lesssim 1$ is satisfied, and,
{as predicted by our theoretical analysis,  we observe a linear convergence in the $H^1$ norm. Moreover, second order convergence in the $L^2$ norm is also observed.}
\begin{table}
\centering
\caption{{Example~\ref{sec:scherk}: computed errors and estimated convergence rates (uniform Voronoi meshes).}}
\label{tab:scherk-uniformVoro}
\begin{tabular}{
c
S[table-format=1.{\roundPrecision}e-1]
S[table-format=6.0]
S[table-format=2.0]
S[table-format=1.{\roundPrecision}e-1]
S[table-format=1.{\roundPrecision}]
S[table-format=1.{\roundPrecision}e-1]
S[table-format=1.{\roundPrecision}]
S[table-format=1.{\roundPrecision}]
S[table-format=1.{\roundPrecision}]
}
\toprule
{Mesh} & {$h$} & {$N$} & {It} & {$e_{H^1}$} & {ecr} & {$e_{L^2}$} & {ecr} & {$C_1$} & {$C_2$}\\
\midrule
u-scherk$_{1}$   &   5.845703e-02   &   4079   &   25   &   2.830468e-02   &   {-}   &   7.509648e-04   &   {-}   &   0.550209   &   2.054193\\
u-scherk$_{2}$   &   3.876636e-02   &   8158   &   27   &   1.972075e-02   &   1.042652   &   3.571190e-04   &   2.144681   &   0.563218   &   1.972077\\
u-scherk$_{3}$   &   2.792674e-02   &   16323   &   29   &   1.396771e-02   &   0.994622   &   1.802687e-04   &   1.971292   &   0.547937   &   1.955020\\
u-scherk$_{4}$   &   1.959732e-02   &   32664   &   31   &   9.785263e-03   &   1.026011   &   8.673115e-05   &   2.109376   &   0.541921   &   1.919408\\
u-scherk$_{5}$   &   1.382014e-02   &   65275   &   29   &   6.986006e-03   &   0.973425   &   4.502900e-05   &   1.893609   &   0.540528   &   1.908553\\
u-scherk$_{6}$   &   1.008418e-02   &   130555   &   30   &   4.922164e-03   &   1.010295   &   2.266426e-05   &   1.980761   &   0.518631   &   1.889713\\
u-scherk$_{7}$   &   7.062998e-03   &   261164   &   31   &   3.484116e-03   &   0.996703   &   1.123266e-05   &   2.024835   &   0.522650   &   1.886500\\
u-scherk$_{8}$   &   5.043997e-03   &   522210   &   30   &   2.472460e-03   &   0.990013   &   5.611956e-06   &   2.002901   &   0.517134   &   1.884952\\
\bottomrule
\end{tabular}
\end{table}
\begin{table}
\centering
\caption{{Example~\ref{sec:scherk}: computed errors and estimated convergence rates (random Voronoi meshes).}}
\label{tab:scherk-randVoro}
\begin{tabular}{
c
S[table-format=1.{\roundPrecision}e-1]
S[table-format=6.0]
S[table-format=2.0]
S[table-format=1.{\roundPrecision}e-1]
S[table-format=1.{\roundPrecision}]
S[table-format=1.{\roundPrecision}e-1]
S[table-format=1.{\roundPrecision}]
S[table-format=1.{\roundPrecision}]
S[table-format=1.{\roundPrecision}]
}
\toprule
{Mesh} & {$h$} & {$N$} & {It} & {$e_{H^1}$} & {ecr} & {$e_{L^2}$} & {ecr} & {$C_1$} & {$C_2$}\\
\midrule
scherk$_{1}$   &   1.182205e-01   &   4098   &   29   &   3.544840e-02   &   {-}   &   1.450268e-03   &   {-}   &   0.312319   &   2.363619\\
scherk$_{2}$   &   7.672309e-02   &   8197   &   31   &   2.444451e-02   &   1.072232   &   6.571025e-04   &   2.283857   &   0.333482   &   2.316465\\
scherk$_{3}$   &   5.475856e-02   &   16391   &   34   &   1.761149e-02   &   0.946236   &   3.761298e-04   &   1.610201   &   0.328400   &   2.302282\\
scherk$_{4}$   &   4.149791e-02   &   32778   &   34   &   1.254106e-02   &   0.979887   &   1.727792e-04   &   2.244998   &   0.310422   &   2.332219\\
scherk$_{5}$   &   2.848743e-02   &   65551   &   33   &   8.769878e-03   &   1.032176   &   8.586859e-05   &   2.017676   &   0.316225   &   2.306425\\
scherk$_{6}$   &   2.126006e-02   &   131087   &   36   &   6.221969e-03   &   0.990535   &   4.452327e-05   &   1.895455   &   0.300526   &   2.313271\\
scherk$_{7}$   &   1.751864e-02   &   262167   &   35   &   4.353764e-03   &   1.030256   &   2.098014e-05   &   2.171154   &   0.255559   &   2.292345\\
scherk$_{8}$   &   1.078938e-02   &   524326   &   32   &   3.101653e-03   &   0.978474   &   1.076297e-05   &   1.925938   &   0.293896   &   2.296097\\
\bottomrule
\end{tabular}
\end{table}

\subsection{Test 4}\label{sec:schwarzD}
\begin{figure}
\centering
\subfloat{\includegraphics[width=0.5\textwidth]{figures/example_unitsquare_randVoro}}
\subfloat{\includegraphics[width=0.5\textwidth]{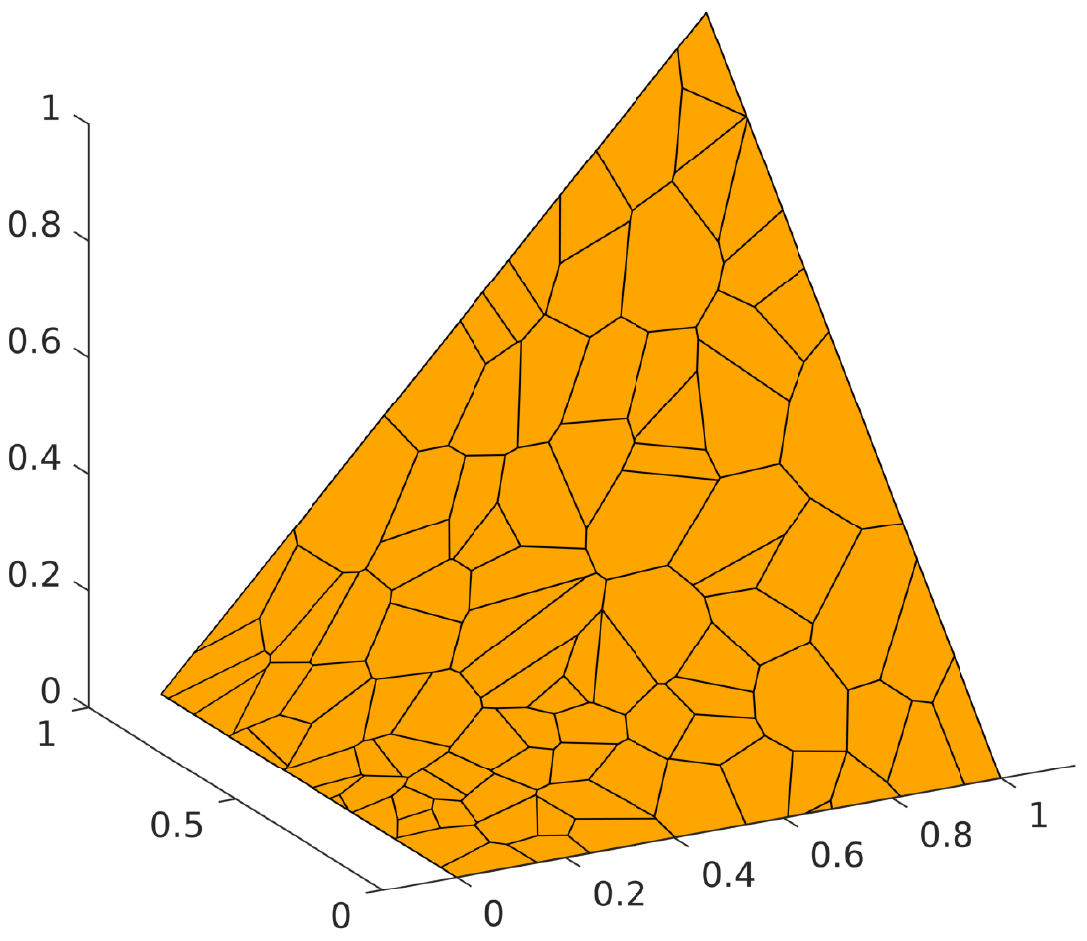}}
\caption{{Example~\ref{sec:schwarzD}: example of the computational mesh (left) and corresponding computed solution (right).}}
\label{fig:example-schwarzD-patch}
\end{figure}

\begin{figure}
\centering
\includegraphics[width=0.618\textwidth]{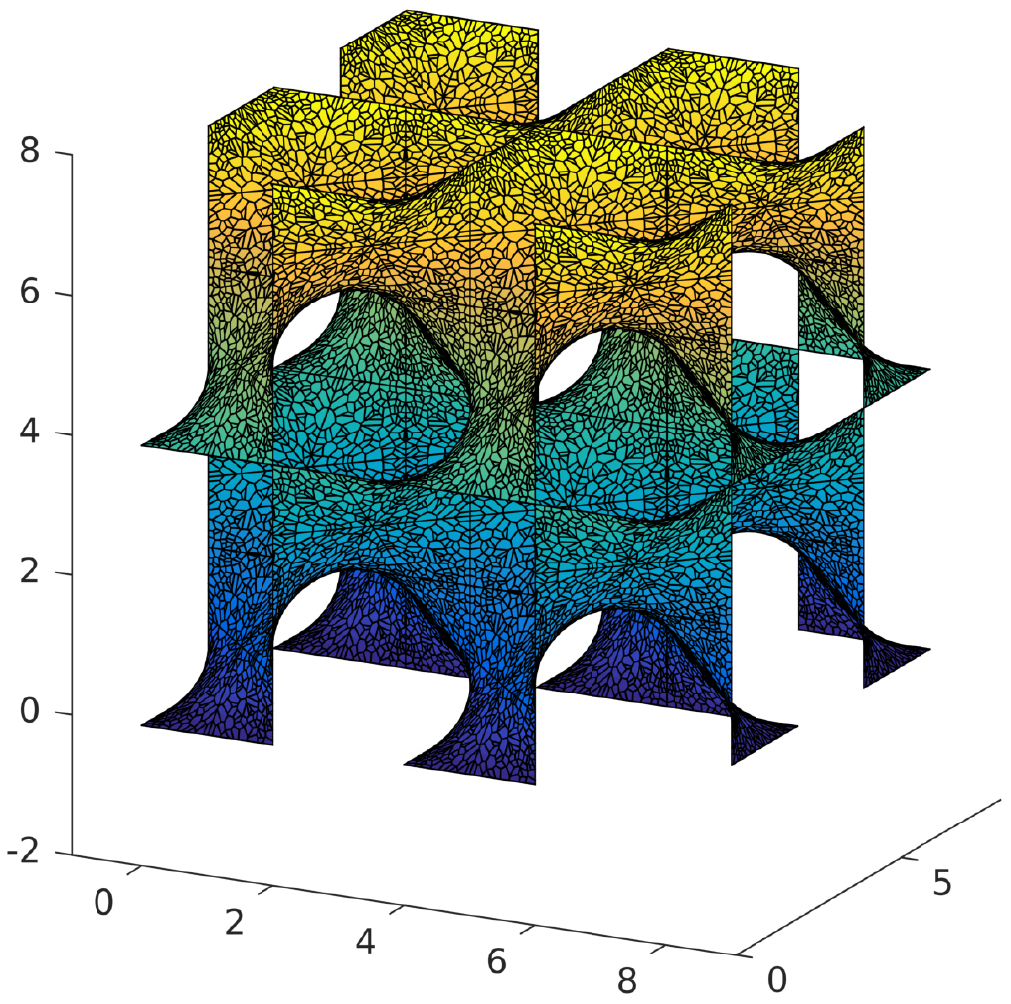}
\caption{{Example~\ref{sec:schwarzD}: representation of the Schwarz D surface obtained by rotating and translating the \emph{patch} shown in Figure~\ref{fig:example-schwarzD-patch}.}}
\label{fig:example-schwarzD}
\end{figure}

The minimal surface problem \eqref{VEM:minpb} is solved on $\Omega=(0,1)^2$ with the following boundary conditions
\[
\begin{cases}
\varphi = 0 & \text{ on }y = 0\text{ and }x = 0,\\
\varphi= x & \text{ on }y = 1,\\
\varphi = y & \text{ on }x = 1.
\end{cases}
\]
A typical solution on a coarse mesh is shown in Figure~\ref{fig:example-schwarzD-patch}. We recall that by properly rotating and translating this minimal surface, it is possible to obtain the so-called Schwarz D surface (see Figure~\ref{fig:example-schwarzD}). Results on uniform and random Voronoi meshes are shown in Tables~\ref{tab:schwarz-uniformVoro} and \ref{tab:schwarz-randVoro}, respectively. The reference FEM solution is computed on a Delaunay triangular mesh with $7767583$ nodes and $15524627$ triangles. 
{Also in this case we observe $\Cuh \lesssim 1$, a linear convergence in the $H^1$ norm, and a quadratic convergence in the $L^2$ norm}.
\begin{table}
\centering
\caption{{Example~\ref{sec:schwarzD}: computed errors and estimated convergence rates (uniform Voronoi meshes).}}
\label{tab:schwarz-uniformVoro}
\begin{tabular}{
c
S[table-format=1.{\roundPrecision}e-1]
S[table-format=6.0]
S[table-format=2.0]
S[table-format=1.{\roundPrecision}e-1]
S[table-format=1.{\roundPrecision}]
S[table-format=1.{\roundPrecision}e-1]
S[table-format=1.{\roundPrecision}]
S[table-format=1.{\roundPrecision}]
S[table-format=1.{\roundPrecision}]
}
\toprule
{Mesh} & {$h$} & {$N$} & {It} & {$e_{H^1}$} & {ecr} & {$e_{L^2}$} & {ecr} & {$C_1$} & {$C_2$}\\
\midrule
u-square$_{1}$   &   3.400609e-02   &   4077   &   22   &   1.138459e-02   &   {-}   &   1.118740e-04   &   {-}   &   0.229888   &   0.499163\\
u-square$_{2}$   &   2.698981e-02   &   8155   &   22   &   8.004869e-03   &   1.016086   &   5.503866e-05   &   2.046354   &   0.203007   &   0.494791\\
u-square$_{3}$   &   1.740063e-02   &   16325   &   22   &   5.697379e-03   &   0.979859   &   2.781825e-05   &   1.966220   &   0.222403   &   0.494462\\
u-square$_{4}$   &   1.211312e-02   &   32636   &   22   &   4.022932e-03   &   1.004724   &   1.387913e-05   &   2.007472   &   0.225823   &   0.494165\\
u-square$_{5}$   &   8.695573e-03   &   65297   &   22   &   2.848237e-03   &   0.995806   &   6.894715e-06   &   2.017594   &   0.221170   &   0.491442\\
u-square$_{6}$   &   6.110583e-03   &   130532   &   22   &   2.021494e-03   &   0.989973   &   3.443489e-06   &   2.004614   &   0.221452   &   0.488900\\
u-square$_{7}$   &   4.472890e-03   &   261135   &   22   &   1.437691e-03   &   0.982950   &   1.721962e-06   &   1.998851   &   0.213786   &   0.488652\\
u-square$_{8}$   &   3.113299e-03   &   522236   &   22   &   1.028641e-03   &   0.966120   &   8.599788e-07   &   2.003549   &   0.217132   &   0.488514\\
\bottomrule
\end{tabular}
\end{table}
\begin{table}
\centering
\caption{{Example~\ref{sec:schwarzD}: computed errors and estimated convergence rates (random Voronoi meshes).}}
\label{tab:schwarz-randVoro}
\begin{tabular}{
c
S[table-format=1.{\roundPrecision}e-1]
S[table-format=6.0]
S[table-format=2.0]
S[table-format=1.{\roundPrecision}e-1]
S[table-format=1.{\roundPrecision}]
S[table-format=1.{\roundPrecision}e-1]
S[table-format=1.{\roundPrecision}]
S[table-format=1.{\roundPrecision}]
S[table-format=1.{\roundPrecision}]
}
\toprule
{Mesh} & {$h$} & {$N$} & {It} & {$e_{H^1}$} & {ecr} & {$e_{L^2}$} & {ecr} & {$C_1$} & {$C_2$}\\
\midrule
square$_{1}$   &   6.384666e-02   &   5006   &   22   &   1.282674e-02   &   {-}   &   1.712198e-04   &   {-}   &   0.135763   &   0.613291\\
square$_{2}$   &   4.344562e-02   &   10008   &   22   &   9.094385e-03   &   0.992785   &   8.512310e-05   &   2.017618   &   0.142351   &   0.618701\\
square$_{3}$   &   3.470002e-02   &   20007   &   22   &   6.384179e-03   &   1.021613   &   4.187674e-05   &   2.048132   &   0.125849   &   0.617691\\
square$_{4}$   &   2.405393e-02   &   40011   &   22   &   4.511768e-03   &   1.001724   &   2.082425e-05   &   2.015988   &   0.127508   &   0.613496\\
square$_{5}$   &   1.726980e-02   &   80007   &   22   &   3.210831e-03   &   0.981758   &   1.058687e-05   &   1.952504   &   0.126651   &   0.618670\\
square$_{6}$   &   1.140086e-02   &   160028   &   22   &   2.268962e-03   &   1.001702   &   5.246371e-06   &   2.025514   &   0.135102   &   0.616167\\
square$_{7}$   &   8.860795e-03   &   320020   &   22   &   1.613746e-03   &   0.983398   &   2.630185e-06   &   1.992635   &   0.123153   &   0.617315\\
square$_{8}$   &   6.248831e-03   &   640035   &   22   &   1.149771e-03   &   0.978144   &   1.308687e-06   &   2.014110   &   0.123295   &   0.616375\\
\bottomrule
\end{tabular}
\end{table}


\subsection{Test 5}\label{sec:matlab}
\begin{figure}
\centering
\subfloat{\includegraphics[width=0.5\textwidth]{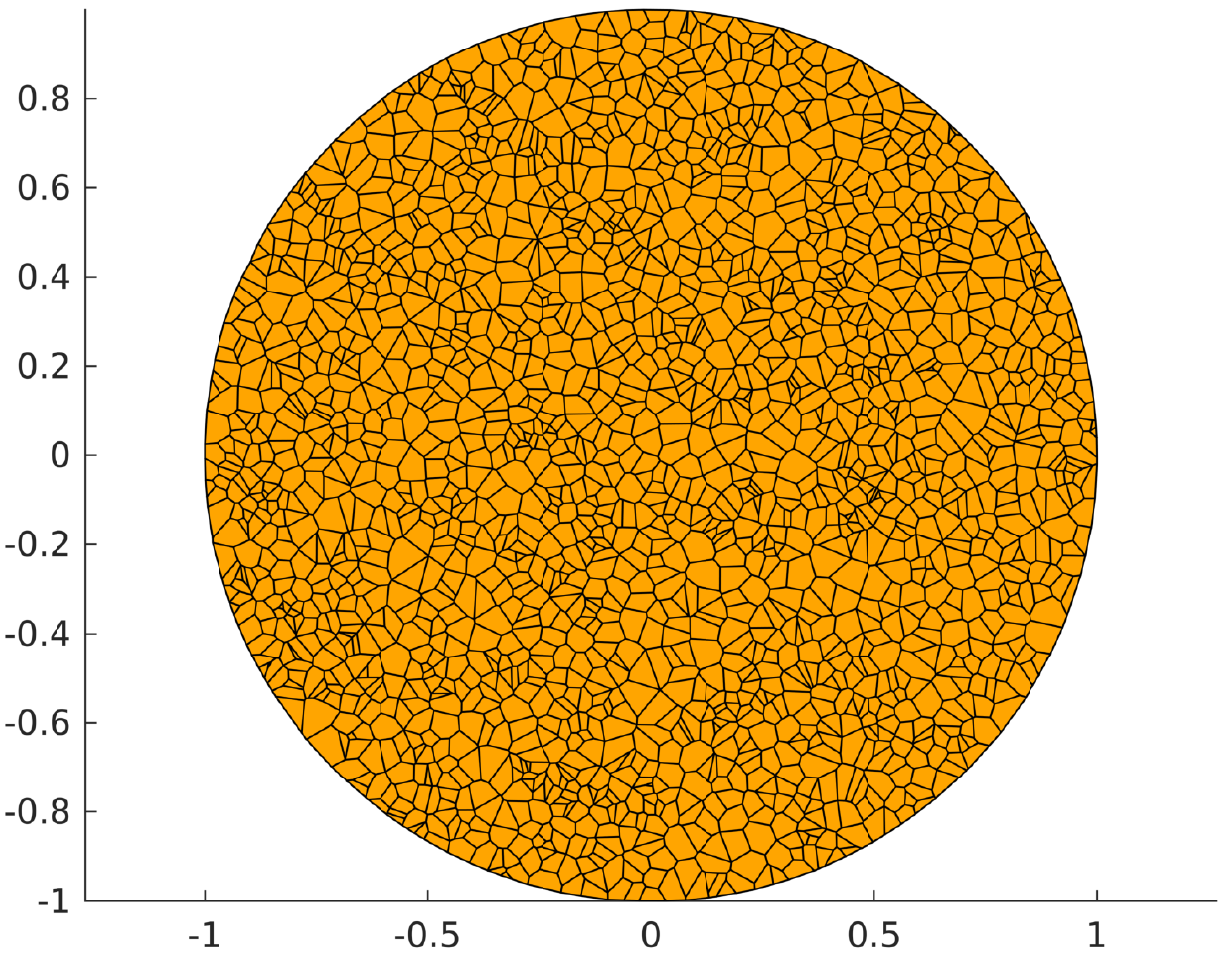}}
\subfloat{\includegraphics[width=0.5\textwidth]{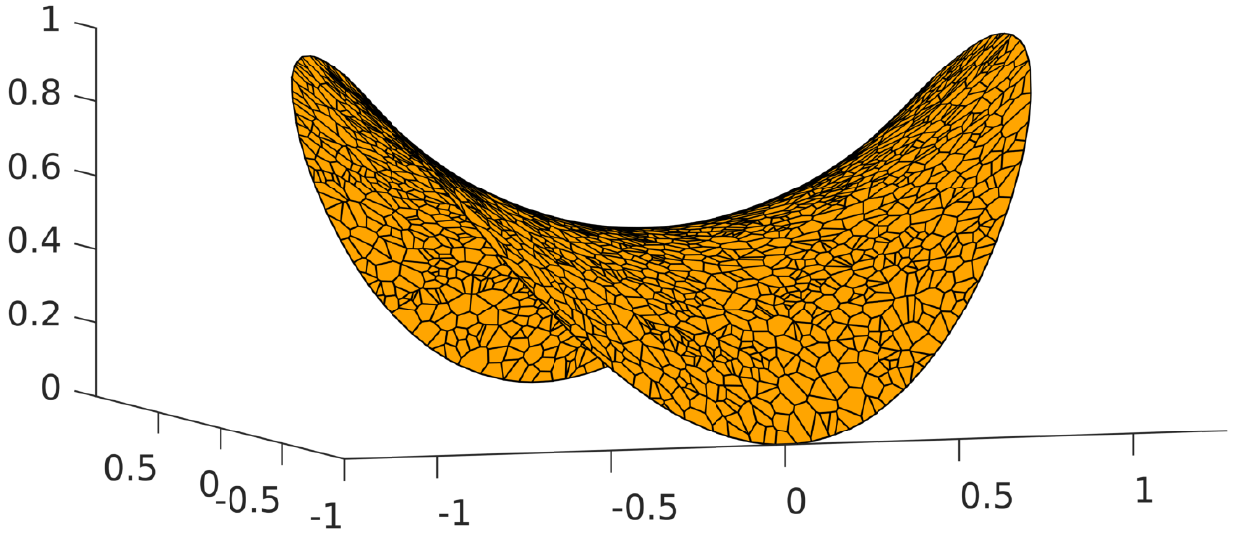}}
\caption{{Example~\ref{sec:matlab}: example of the computational mesh (left) and corresponding computed solution (right).}}
\label{fig:example-matlab}
\end{figure}

Here we consider a minimal surface problem on the unit disk, where the boundary condition is $\varphi(x,y) = x^2$. 
A typical solution on a coarse mesh is shown in Figure~\ref{fig:example-matlab}. Results on uniform and random
Voronoi meshes are shown in Tables~\ref{tab:matlab-uniformVoro}
and \ref{tab:matlab-randVoro}, respectively. 
{Again, $\Cuh \lesssim 1$ and we observe a linear convergence in the $H^1$ norm. Moreover, second order convergence in the $L^2$ norm is also observed.}
\begin{table}
\centering
\caption{{Example~\ref{sec:matlab}: computed errors and estimated convergence rates (uniform Voronoi meshes).}}
\label{tab:matlab-uniformVoro}
\begin{tabular}{
c
S[table-format=1.{\roundPrecision}e-1]
S[table-format=6.0]
S[table-format=2.0]
S[table-format=1.{\roundPrecision}e-1]
S[table-format=1.{\roundPrecision}]
S[table-format=1.{\roundPrecision}e-1]
S[table-format=1.{\roundPrecision}]
S[table-format=1.{\roundPrecision}]
S[table-format=1.{\roundPrecision}]
}
\toprule
{Mesh} & {$h$} & {$N$} & {It} & {$e_{H^1}$} & {ecr} & {$e_{L^2}$} & {ecr} & {$C_1$} & {$C_2$}\\
\midrule
u-circ$_{1}$   &   6.167744e-02   &   4074   &   22   &   2.491744e-02   &   {-}   &   2.175920e-04   &   {-}   &   0.416637   &   1.640193\\
u-circ$_{2}$   &   4.246809e-02   &   8166   &   22   &   1.759650e-02   &   1.000548   &   1.086885e-04   &   1.996496   &   0.427784   &   1.641695\\
u-circ$_{3}$   &   3.054830e-02   &   16303   &   22   &   1.243710e-02   &   1.003850   &   5.396769e-05   &   2.025256   &   0.422185   &   1.646733\\
u-circ$_{4}$   &   2.163149e-02   &   32619   &   22   &   8.797412e-03   &   0.998425   &   2.688975e-05   &   2.008924   &   0.416159   &   1.625852\\
u-circ$_{5}$   &   1.506714e-02   &   65275   &   23   &   6.232316e-03   &   0.993809   &   1.346466e-05   &   1.994126   &   0.425050   &   1.636229\\
u-circ$_{6}$   &   1.064505e-02   &   130546   &   23   &   4.427231e-03   &   0.986771   &   6.718177e-06   &   2.006161   &   0.425504   &   1.636564\\
u-circ$_{7}$   &   7.731530e-03   &   261069   &   23   &   3.158871e-03   &   0.974115   &   3.354122e-06   &   2.004524   &   0.414419   &   1.637127\\
u-circ$_{8}$   &   5.459296e-03   &   522188   &   23   &   2.274585e-03   &   0.947480   &   1.675214e-06   &   2.002903   &   0.414589   &   1.635565\\
\bottomrule
\end{tabular}
\end{table}
\begin{table}
\centering
\caption{{Example~\ref{sec:matlab}: computed errors and estimated convergence rates (random Voronoi meshes).}}
\label{tab:matlab-randVoro}
\begin{tabular}{
c
S[table-format=1.{\roundPrecision}e-1]
S[table-format=6.0]
S[table-format=2.0]
S[table-format=1.{\roundPrecision}e-1]
S[table-format=1.{\roundPrecision}]
S[table-format=1.{\roundPrecision}e-1]
S[table-format=1.{\roundPrecision}]
S[table-format=1.{\roundPrecision}]
S[table-format=1.{\roundPrecision}]
}
\toprule
{Mesh} & {$h$} & {$N$} & {It} & {$e_{H^1}$} & {ecr} & {$e_{L^2}$} & {ecr} & {$C_1$} & {$C_2$}\\
\midrule
circ$_{1}$   &   1.207272e-01   &   4238   &   21   &   3.111517e-02   &   {-}   &   4.282237e-04   &   {-}   &   0.271203   &   2.131472\\
circ$_{2}$   &   8.719227e-02   &   8398   &   22   &   2.187564e-02   &   1.030328   &   2.050550e-04   &   2.153430   &   0.264488   &   2.113351\\
circ$_{3}$   &   6.233188e-02   &   16662   &   22   &   1.545690e-02   &   1.013865   &   1.021701e-04   &   2.033575   &   0.260636   &   2.097044\\
circ$_{4}$   &   4.554244e-02   &   33175   &   22   &   1.094648e-02   &   1.002048   &   5.069568e-05   &   2.035236   &   0.253473   &   2.102582\\
circ$_{5}$   &   3.382851e-02   &   66096   &   23   &   7.744545e-03   &   1.003984   &   2.544399e-05   &   2.000144   &   0.239361   &   2.081728\\
circ$_{6}$   &   2.421469e-02   &   131883   &   23   &   5.479975e-03   &   1.001403   &   1.264160e-05   &   2.025128   &   0.237022   &   2.084310\\
circ$_{7}$   &   1.804084e-02   &   263268   &   23   &   3.898897e-03   &   0.984892   &   6.310083e-06   &   2.010378   &   0.224970   &   2.082481\\
circ$_{8}$   &   1.199765e-02   &   525901   &   23   &   2.783893e-03   &   0.973619   &   3.121175e-06   &   2.034681   &   0.238740   &   2.077179\\
\bottomrule
\end{tabular}
\end{table}


\subsection{Test 6}\label{sec:cantor}
\begin{figure}
\centering
\subfloat{\includegraphics[width=0.5\textwidth]{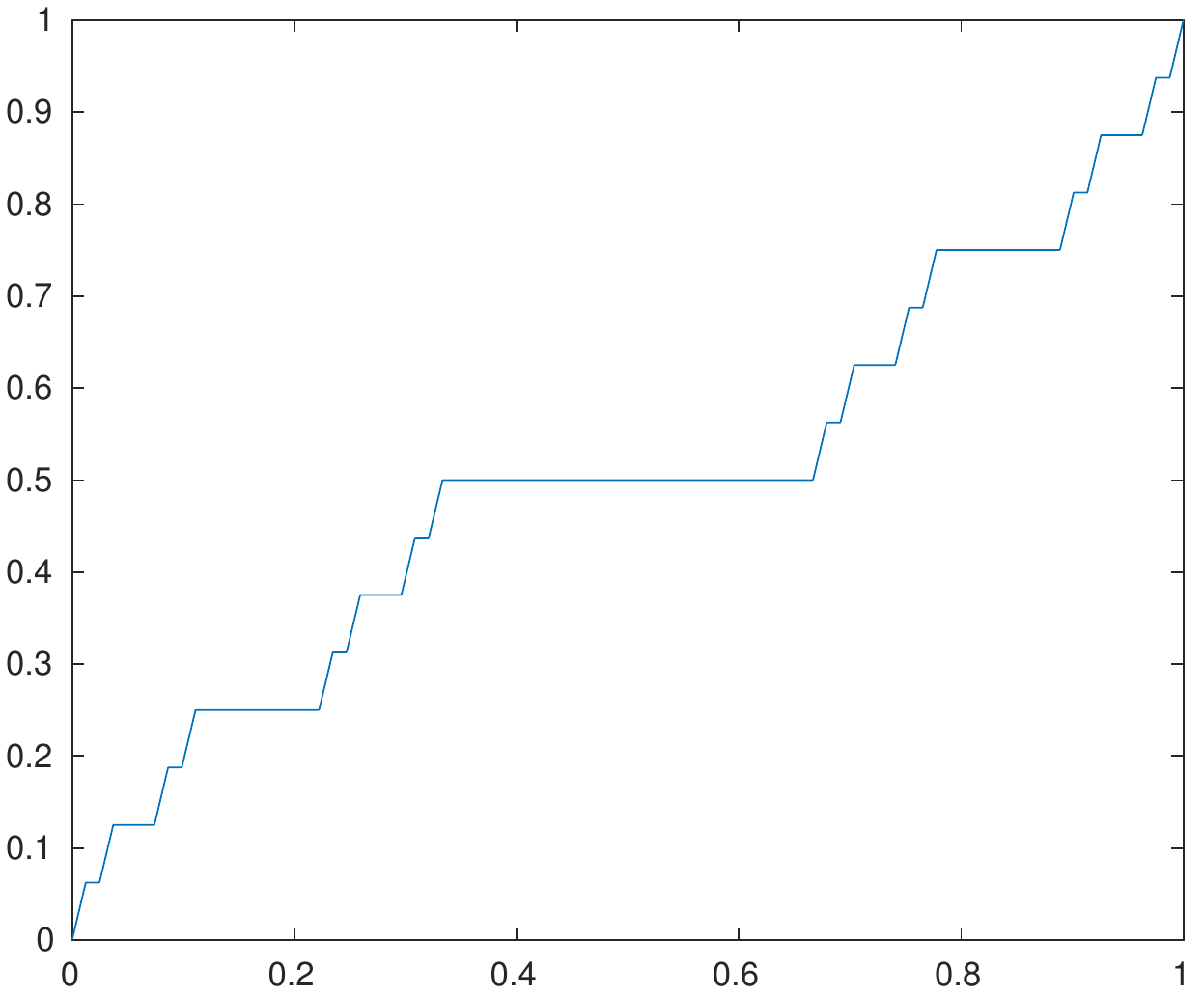}}
\subfloat{\includegraphics[width=0.5\textwidth]{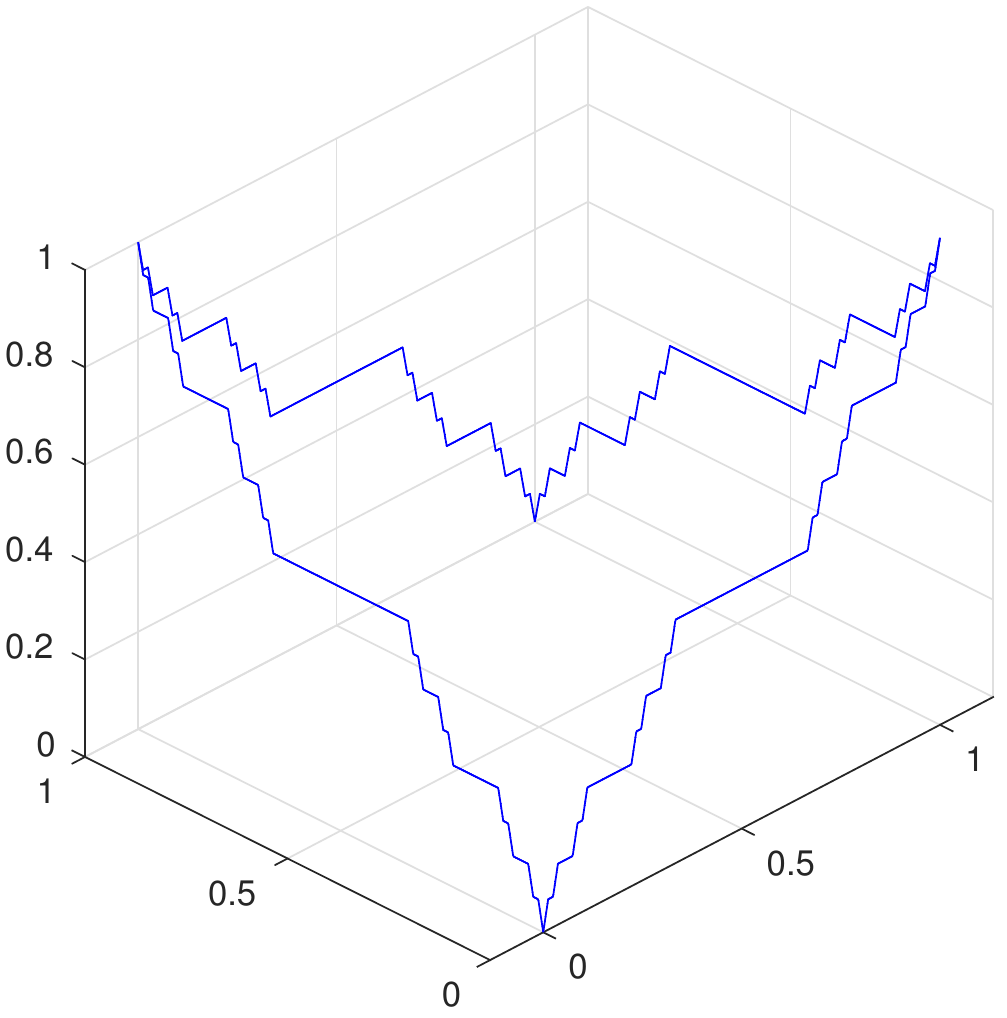}}
\caption{{Example~\ref{sec:cantor}: fourth iterate of a sequence of function converging to the Cantor function on $[0,1]$ (left) which is used as as boundary condition for Example~\ref{sec:cantor} (right).}}
\label{fig:cantor-bc}
\end{figure}

\begin{figure}
\centering
\subfloat{\includegraphics[width=0.5\textwidth]{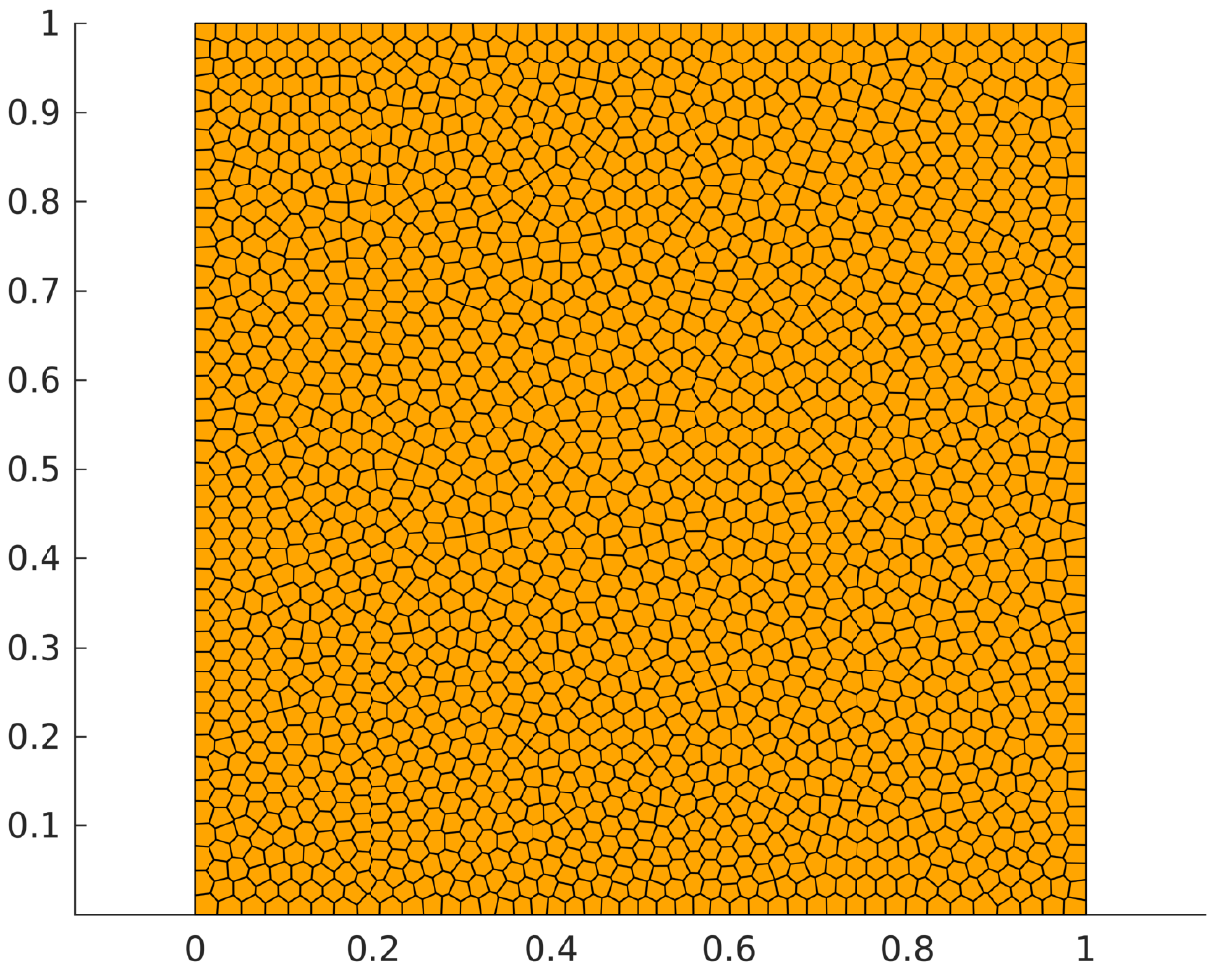}}
\subfloat{\includegraphics[width=0.5\textwidth]{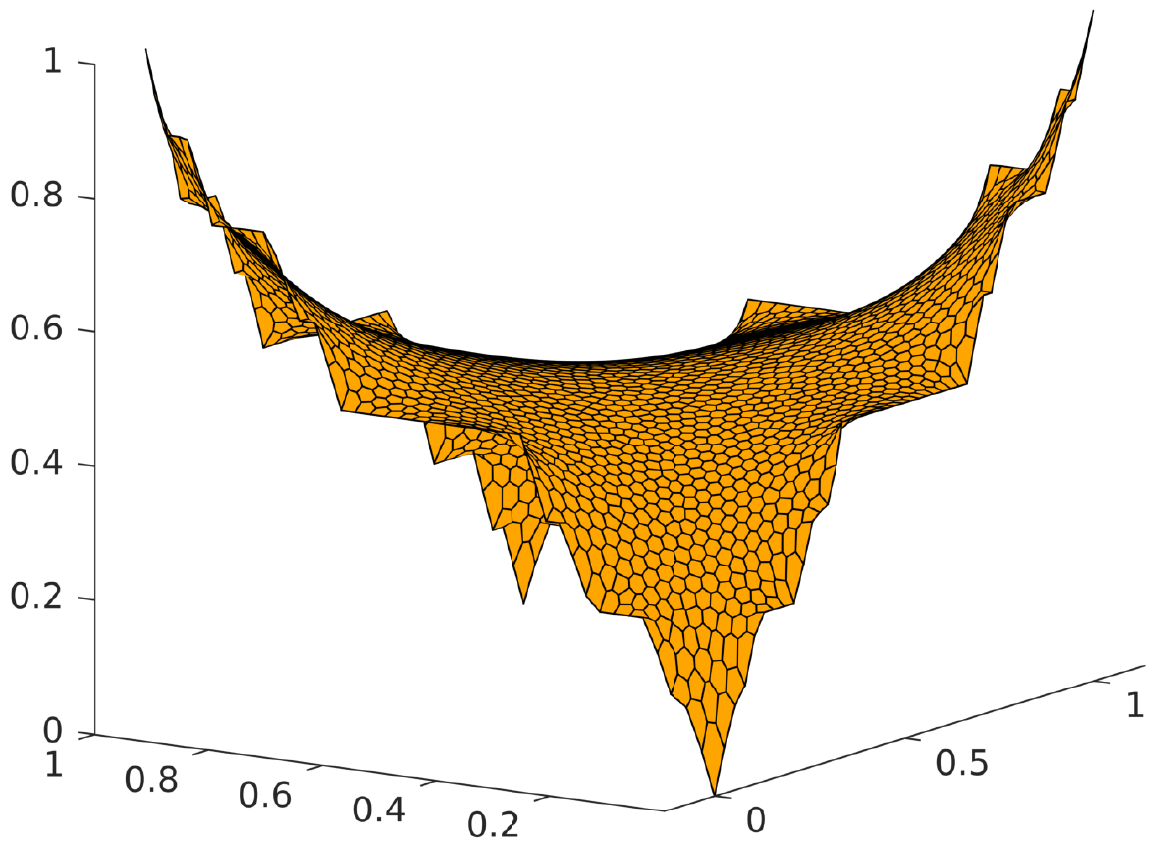}}
\caption{{Example~\ref{sec:cantor}: example of the computational mesh (left) and corresponding computed solution (right).}}
\label{fig:example-cantor}
\end{figure}

In the last example, the minimal surface problem is again solved on $\Omega=(0,1)^2$. As Dirichlet boundary conditions, we require the solution to match proper reflections of the fourth iterate of a sequence of functions converging to the Cantor function (see Figure~\ref{fig:cantor-bc}). Note that the exact solution does not satisfy the regularity assumptions of Theorem \ref{th:main}.
A typical solution on a coarse mesh is shown in Figure~\ref{fig:example-cantor}. Results on uniform and random Voronoi meshes are shown in Tables~\ref{tab:cantor-uniformVoro} and \ref{tab:cantor-randVoro}, respectively. The reference FEM solution is computed on a Delaunay triangular mesh with $7768041$ nodes and $15525051$ triangles. Such mesh is constructed in order to have all the nodes where the Dirichlet data is just continuous as boundary nodes. Observe that the assumption $\Cuh \lesssim 1$ does not hold in this case. This example shows that a lack of regularity in the boundary data may severely affect the convergence properties of the method.
\begin{table}
\centering
\caption{{Example~\ref{sec:cantor}: computed errors and estimated convergence rates (uniform Voronoi meshes).}}
\label{tab:cantor-uniformVoro}
\begin{tabular}{
c
S[table-format=1.{\roundPrecision}e-1]
S[table-format=6.0]
S[table-format=3.0]
S[table-format=1.{\roundPrecision}e-1]
S[table-format=1.{\roundPrecision}]
S[table-format=1.{\roundPrecision}e-1]
S[table-format=1.{\roundPrecision}]
S[table-format=2.{\roundPrecision}]
S[table-format=2.{\roundPrecision}]
}
\toprule
{Mesh} & {$h$} & {$N$} & {It} & {$e_{H^1}$} & {ecr} & {$e_{L^2}$} & {ecr} & {$C_1$} & {$C_2$}\\
\midrule
u-square$_{1}$   &   3.400609e-02   &   4077   &   31   &   4.609498e-01   &   {-}   &   3.013669e-02   &   {-}   &   3.374811   &   7.327850\\
u-square$_{2}$   &   2.698981e-02   &   8155   &   45   &   4.194200e-01   &   0.272380   &   2.415286e-02   &   0.638541   &   3.983681   &   9.709486\\
u-square$_{3}$   &   1.740063e-02   &   16325   &   45   &   3.902435e-01   &   0.207767   &   2.114981e-02   &   0.382591   &   6.553732   &   14.570698\\
u-square$_{4}$   &   1.211312e-02   &   32636   &   46   &   2.826745e-01   &   0.931042   &   4.374640e-03   &   4.549632   &   7.861173   &   17.202510\\
u-square$_{5}$   &   8.695573e-03   &   65297   &   59   &   3.419749e-01   &   -0.549193   &   1.280841e-02   &   -3.097998   &   9.116705   &   20.257355\\
u-square$_{6}$   &   6.110583e-03   &   130532   &   67   &   3.683411e-01   &   -0.214451   &   1.364939e-02   &   -0.183617   &   10.450648   &   23.071953\\
u-square$_{7}$   &   4.472890e-03   &   261135   &   84   &   2.695212e-01   &   0.900934   &   6.281133e-03   &   2.238602   &   11.192523   &   25.582839\\
u-square$_{8}$   &   3.113299e-03   &   522236   &   100   &   3.079599e-01   &   -0.384724   &   7.085787e-03   &   -0.347839   &   12.481802   &   28.082226\\
\bottomrule
\end{tabular}
\end{table}

\begin{table}
\centering
\caption{{Example~\ref{sec:cantor}: computed errors and estimated convergence rates (random Voronoi meshes).}}
\label{tab:cantor-randVoro}
\begin{tabular}{
c
S[table-format=1.{\roundPrecision}e-1]
S[table-format=6.0]
S[table-format=3.0]
S[table-format=1.{\roundPrecision}e-1]
S[table-format=1.{\roundPrecision}]
S[table-format=1.{\roundPrecision}e-1]
S[table-format=1.{\roundPrecision}]
S[table-format=2.{\roundPrecision}]
S[table-format=2.{\roundPrecision}]
}
\toprule
{Mesh} & {$h$} & {$N$} & {It} & {$e_{H^1}$} & {ecr} & {$e_{L^2}$} & {ecr} & {$C_1$} & {$C_2$}\\
\midrule
square$_{1}$   &   6.384666e-02   &   5006   &   44   &   5.381897e-01   &   {-}   &   3.317278e-02   &   {-}   &   2.538615   &   11.467810\\
square$_{2}$   &   4.344562e-02   &   10008   &   50   &   5.372790e-01   &   0.004889   &   3.101866e-02   &   0.193839   &   4.023503   &   17.487349\\
square$_{3}$   &   3.470002e-02   &   20007   &   54   &   5.743910e-01   &   -0.192849   &   5.022859e-02   &   -1.391648   &   3.812010   &   18.710040\\
square$_{4}$   &   2.405393e-02   &   40011   &   57   &   3.932863e-01   &   1.093025   &   8.672985e-03   &   5.068365   &   5.212454   &   25.079448\\
square$_{5}$   &   1.726980e-02   &   80007   &   70   &   4.376071e-01   &   -0.308196   &   1.380389e-02   &   -1.341311   &   6.139793   &   29.991971\\
square$_{6}$   &   1.140086e-02   &   160028   &   85   &   4.525599e-01   &   -0.096933   &   2.115745e-02   &   -1.232026   &   8.809973   &   40.180009\\
square$_{7}$   &   8.860795e-03   &   320020   &   97   &   3.127613e-01   &   1.066267   &   7.969529e-03   &   2.817655   &   8.732460   &   43.772147\\
square$_{8}$   &   6.248831e-03   &   640035   &   113   &   3.070381e-01   &   0.053289   &   6.110790e-03   &   0.766280   &   10.181842   &   50.901080\\
\bottomrule
\end{tabular}
\end{table}


\section{Conclusions}
We presented the lowest order {V}irtual {E}lement discretization of a minimal surface problem. Optimal error estimate in the $H^1$-norm has been derived and several numerical tests assessing the validity of the theoretical results have been presented. Moreover, the convergence properties in the $L^2$-norm has been numerically investigated.
\section*{Acknowledgments}
The authors are members of the INdAM Research group GNCS and this work is partially funded by INDAM-GNCS. P.F.A. and M.V. acknowledge the financial support of MIUR thourgh the PRIN grant n. 201744KLJL. 

\bibliographystyle{abbrv}
\bibliography{biblio}

\end{document}